\documentclass[11pt]{amsart}
\RequirePackage{newlfont}
\usepackage{amscd}
\usepackage{amsmath}
\usepackage{amsfonts}
\usepackage{amssymb}
\usepackage{enumerate}
\usepackage{graphicx}
\usepackage{latexsym}
\usepackage{color}
\usepackage{bbm}

\newtheorem{theorem}{Theorem}[section]
\newtheorem{lemma}[theorem]{Lemma}
\newtheorem{proposition}[theorem]{Proposition}
\newtheorem{corollary}[theorem]{Corollary}
\theoremstyle{definition}
\newtheorem{definition}[theorem]{Definition}
\newtheorem{example}[theorem]{Example}

\newtheorem{remark}[theorem]{Remark}

\theoremstyle{remark}

\def\supp{\textup{supp }}
\def\e{\textup{e}}

\numberwithin{equation}{section}

\begin{document}

	\bibliographystyle{plainnat}

	\setcounter{page}{1}
	
	\title[Dirichlet problem]{The non-homogeneous Dirichlet problem for the $p(x)$-Laplacian with unbounded $p(x)$ on a smooth domain.}
	\author[M. A. Khamsi, J. Lang, O. M\'{e}ndez \& A. Nekvinda]{Mohamed  A. Khamsi, Jan Lang,  Osvaldo M\'{e}ndez, Ale\v{s} Nekvinda}

	\address{Mohamed A. Khamsi\\ Department of Applied Mathematics and Sciences, Khalifa University, Abu Dhabi, UAE}
	\email{mohamed.khamsi@ku.ac.ae}
	\address{Jan Lang\\ Department of Mathematics, 100 Math Tower, 231 West 18th Ave., Columbus,
		OH 43210-1174, USA \\ and \\ Department of Mathematics, Faculty of Electrical Engineering, Czech Technical University in Prague,   Technick\'a~2, 166~27 Praha~6, Czech Republic}
	\email{lang.162@osu.edu}
	\address{Osvaldo M\'{e}ndez\\Department of Mathematical Sciences, The University of Texas at El Paso, El Paso, TX 79968, USA}
	\email{osmendez@utep.edu}
	\address{Ale\v{s} Nekvind\\Faculty of Civil Engineering, Czech Technical University, Thakurova 7, 166 29 Prague 6, Czech Republic}
	\email{ales.nekvinda@cvut.cz}
	
	
	\thanks{The fourth author was supported by the grant P202/23-04720S of the Grant Agency of the Czech Republic}

	\begin{abstract}
	
	This paper examines the solvability of the Dirichlet problem for the variable exponent \( p \)-Laplacian in the case of unbounded \( p(x) \). For a bounded domain \( \Omega \subset \mathbb{R}^n \) with a smooth boundary and \( p \) satisfying \( p_- = \text{essinf}_{\Omega} p > n \) and $\varphi$ in the Sobolev space $W^{1,p(x)}(\Omega)$, we investigate the problem
	\[
	\Delta_p(u) = 0 \ \text{in} \ \Omega, \quad u|_{\partial \Omega} = \varphi.
	\]
	We introduce the space \( V_0^{1,p}(\Omega) \), which is the natural solution space for the minimization of the Dirichlet integral given the unbounded nature of \( p(x) \). Our main results establish the existence and uniqueness of solutions within this space. Since $V^{1,p}_0(\Omega)$ is not defined via a TVS topology, the paper includes the description of the necessary modular topological framework and discusses Clarkson-type inequalities for unbounded variable exponents, which are interesting in their own right.
	\end{abstract}
	
	\maketitle
	\section{Introduction}
This work is devoted to the analysis of the solvability of the Dirichlet problem for the variable exponent \( p \)-Laplacian in the novel context of an unbounded exponent \( p(x) \). Specifically, let \( \Omega \subset \mathbb{R}^n \) be a bounded domain with a smooth boundary \( \partial \Omega \) ( a \( C^2 \) condition on the boundary suffices), and let \( p \) be a measurable function on \( \Omega \) with \( p_- = \text{essinf}_{\Omega} p > n \) and \( \varphi \in W^{1,p}(\Omega) \). We investigate the solvability of the problem
\begin{equation}\label{DP}
	\begin{cases}
		\Delta_p(u) = \text{div}\left( |\nabla u|^{p-2} \nabla u \right) = 0 & \text{in} \ \Omega, \\
		u|_{\partial \Omega} = \varphi
	\end{cases}
\end{equation}
in \( W^{1,p}(\Omega) \).

For bounded exponents \( p \), the boundary condition in (\ref{DP}) is typically interpreted as \( u - \varphi \in W^{1,p}_0(\Omega) \), where \( W^{1,p}_0(\Omega) \) is the closure of \( C^{\infty}_0(\Omega) \) in the \( W^{1,p}(\Omega) \)-norm. The solvability of (\ref{DP}) in this context follows from the uniform convexity of the Sobolev norm, achieved via the minimization of the Dirichlet integral
\begin{equation}\label{DI}
	D(u) = \int_{\Omega} |\nabla (u - \varphi)|^p \, dx
\end{equation}
on the subspace \( W^{1,p}_0(\Omega) \), where \( u \) is the minimizer of (\ref{DI}).

In sharp contrast, the Sobolev norm is not uniformly convex when the exponent \( p(x) \) is unbounded \cite{Luk}. Addressing the Dirichlet problem (\ref{DP}) in this scenario necessitates a different approach. We demonstrate that the absence of uniform convexity in the Sobolev norm does not preclude the good  geometric behavior of the Dirichlet integral. Our analysis reveals that the classical homogeneous Sobolev space \( W^{1,p}_0(\Omega) \) is inadequate for the unbounded \( p(x) \) case. When \( p(x) \) is unbounded, the appropriate space for minimizing the Dirichlet integral, which consequently defines the boundary condition, strictly contains \( W^{1,p}_0(\Omega) \). This space is denoted by \( V_0^{1,p}(\Omega) \) (see Definition \ref{defv1p0}). We establish that for unbounded \( p(x) \) on \( \Omega \), a solution to (\ref{DP}) exists in \( V_0^{1,p}(\Omega) \) and that, under natural conditions, the solution is unique.

The paper is organized as follows: Section \ref{preliminaries} presents some known results. Section \ref{s3} introduces the necessary function spaces. Section \ref{s4} discusses the modular topology. 
In Section \ref{s5}, we prove Clarkson-type inequalities for variable exponents. Section \ref{s6} focuses on the minimization of the Dirichlet integral. Section \ref{s7} concludes with the proof of the solvability of problem (\ref{DP}). Finally, in Section \ref{Sectionremark}, we indicate how the ideas introduced in this work can be applied to other generalizations of the \( p \)-Laplacian for unbounded variable exponents \( p(x) \).

	\section{Preliminaries and known results}\label{preliminaries}
	Boundary value problems for differential operators with variable exponents under the condition of boundedness of the exponent have been extensively studied. The next few results (the list is by no means exhaustive) have served as motivation for the ideas expounded in this work. The homogeneous Dirichlet problem
	\begin{equation}\label{DPf}
		\begin{cases}
			\Delta_p(u)=\text{div}\left(|\nabla u|^{p(x)-2}\nabla u\right)=f(x,u)\,\,\,\text{in}\,\,\, \Omega\\
			u|_{\partial \Omega}= 0,
		\end{cases}
	\end{equation}
	has been studied in \cite{FZ} (see also \cite{CF}) for a suitable right-hand side $f$, under the conditions:
	\[ 1< p_-=\inf\limits_{x\in \Omega}p(x)\leq \sup\limits_{x\in \Omega}p(x)=p_+<\infty. \]
	Notice also that the boundary condition in problem (\ref{DPf}) is zero.  Under the same assumptions, in \cite{CF} the authors tackle the problem
	\begin{equation*}
		\begin{cases}
			-\Delta_p(u)+a(x)|u|^{p(x)-2}u=f(x,u)\,\,\,\text{in}\,\,\, \Omega\\
			u|_{\partial \Omega}= 0,
		\end{cases}
	\end{equation*} whereas the general homogeneous boundary value problem
	\begin{equation*}
		\begin{cases}
			\sum\limits_iD_i\left(a_i(x,u)|D_iu|^{p(x)-2}D_iu\right)=f(x,u)\,\,\,\text{in}\,\,\, \Omega\\
			u|_{\partial \Omega}= 0
		\end{cases}
	\end{equation*}
	was solved in \cite{AS}. A similar problem with a mixed boundary condition was considered in \cite{AC}. \\
	A wide collection of similar results in this spirit can be found in \cite{AS,H,MR1, MR2,SU,YW,Z1, ZZ} and the references therein; all such works rely essentially on the assumptions of zero boundary values and on the boundedness of the variable exponent $p(x)$, i.e., on the condition $p_+=\sup\limits_{x\in \Omega}p(x)<\infty$. In \cite{MRU}, the authors consider a sequence of variable exponents $(p_j)$ defined on a bounded domain $\Omega \subset {\mathbb R}^n$ and subject to $n<\alpha \leq \inf\limits_{\Omega}p_j\leq \sup\limits_{\Omega}p_{j}<\infty$ for each $j$, such that $\lim\limits_{j\rightarrow \infty}p_j=\infty$ uniformly on a bounded domain $\Omega$, and show that the sequence of solutions $(u_j)$ of the problems
	\begin{equation}\label{DPn}
		\begin{cases}
			\Delta_{p_j}(u)=\text{div}\left(|\nabla u|^{p_j(x)-2}\nabla u\right)=0\,\,\,\text{in}\,\,\, \Omega\\
			u|_{\partial \Omega}= f,
		\end{cases}
	\end{equation} 
	converges to the unique viscosity solution of a certain boundary value problem (Theorem 1.1 in  \cite{MRU}).\\
	It is to be emphasized here that in problem (\ref{DPn}), each exponent $p_j$ is assumed to be bounded in $\Omega$, that is, $n<\alpha \leq \inf\limits_{x \in \Omega}\ p_j(x) \leq \sup\limits_{x \in \Omega}\ p_j(x) <\infty$ (condition 2.9 in \cite{MRU}), for each $j$,  and that the boundary value is Lipschitz continuous on $\partial \Omega$ with Lipschitz constant less than or equal to one.\\
	In a related work \cite{MRU1}, the same authors show that for a measurable admissible exponent $p(x)$ defined on a bounded, convex, smooth domain $\Omega$ such that $p(x)=\infty$ on a convex, smooth subdomain $D\subset \Omega$ and $p(x)$ is continuously differentiable on $\Omega$, bounded and greater than the dimension $n$ on $\Omega \setminus {\overline D}$, the problem
	\begin{equation}\label{DPin}
		\begin{cases}
			\Delta_{p}(u)=\text{div}\left(|\nabla u|^{p(x)-2}\nabla u\right)=0\,\,\,\text{in}\,\,\, \Omega\setminus \overline{D}\\
			\Delta_{\infty}u=0 \,\,\,\text{in}\,\,\, D\\
			B(u)\,\,\,\text{in}\,\,\, \partial D \cap \Omega\\
			u|_{\partial \Omega}= f,
		\end{cases}
	\end{equation}
	where $B$ is a certain boundary operator, has a solution. \cite[p. 2582]{MRU1}.\\
	
	In this work we consider a smooth, bounded domain $\Omega\subset {\mathbb R}^n$ and an unbounded,  measurable exponent $p$ satisfying $n< \inf\limits_{x \in \Omega}\ p(x)$.

	\section{Function spaces of variable exponent} \label{s3}
	By ${\mathcal P}(\Omega)$ we denote the family of all (measurable) functions on $\Omega \subset  {\mathbb R}^n,$  $p:\Omega \to (1, \infty)$ with $p_-:=\inf\limits_{x\in \Omega}p(x) >n$ and $p_+:=\sup\limits_{x\in \Omega}p(x)= \infty$.\\
	
	Unless specifically mentioned otherwise it will always be assumed that the variable exponent $p(\cdot)$ (or simply $p$) belongs to  ${\mathcal P}(\Omega)$, that is, no continuity or log-H\"{o}lder continuity will be assumed.  Moreover, in the sequel $\Omega\subset {\mathbb R}^n$ will stand for a bounded and regular domain (i.e., the boundary $\partial \Omega$ will be assumed to be locally $C^2$).\\
	On the set ${\mathcal M}$ of extended-real-valued Borel-measurable functions defined on $\Omega$, consider the functional $\rho_{p}:{\mathcal M}\rightarrow [0,\infty]$ defined by
	\begin{equation}\label{def-rho}
		\rho_{p}(u):=\int_{\Omega}\frac{|u(x)|^{p(x)}}{p(x)}\ dx.
	\end{equation}
	Then $\rho_p$ is a convex modular on ${\mathcal M}$, that is
	\begin{enumerate}\label{modularproperties}
		\item[(1)] $\rho_p(u) = 0$ if and only if $u = 0$,
		\item[(2)] $\rho_p(\alpha u) = \rho(x)$, if $|\alpha| =1$,
		\item[(3)] $\rho_p(\alpha u + (1-\alpha) v )\leq \alpha\rho_p(u) + (1-\alpha)\rho_p(v)$, for any $\alpha \in [0,1]$
		and any $u, v \in {\mathcal M}$.
	\end{enumerate}
	The variable exponent Lebesgue space $L^p(\Omega)$ is defined as \cite{DHHR,KR}
	\begin{equation*}
		L^p(\Omega)=\left\{u\in {\mathcal M}:\rho_p(\lambda u)<\infty\,\,\text{for some}\,\,\lambda>0\right\},
	\end{equation*}
	which  is a Banach space when endowed with the Luxemburg norm
	\begin{equation}
		\|u\|_{\rho}:=\inf\left \{\lambda>0:\rho_{p}\left(\frac{u}{\lambda}\right) \leq 1 \right\}.
	\end{equation}
	Analogously, on the set ${\mathcal V}\subset {\mathcal M}$ consisting of those functions whose distributional derivatives belong to ${\mathcal M}$ one has the convex modular $\rho_{1,p}:{\mathcal V}\rightarrow [0,\infty]$ defined by
	\begin{equation}\label{defrho}
		\rho_{1,p}(u):= \rho_p(u)+\rho_p\left(|\nabla u|\right) = \int_{\Omega} \frac{|\nabla u|^{p(x)}}{p(x)}dx + \int_{\Omega} \frac{|u(x)|^{p(x)}}{p(x)}dx ,
	\end{equation}
	where $|\nabla u|$ stands for the Euclidean norm of the gradient of $u$.
	The space $W^{1,p}(\Omega)$ is defined as the class of real-valued, Borel measurable functions $u \in {\mathcal V}$ , for which there exists $\lambda>0$ such that $\rho_{1,p}(\lambda u)<\infty.$  The Luxemburg norm $\|\cdot\|_{1,p}$ associated to the modular $\rho_{1,p}$ is given by
	\begin{equation} \label{Lux Sob}
		\|u\|_{1,p}=\inf\left \{\lambda>0:\rho_{1,p}\left(\frac{u}{\lambda}\right) \leq 1 \right\},
	\end{equation}
	with respect to which, $W^{1,p}(\Omega)$ is a Banach space.
	
	\begin{remark}{\normalfont The Luxemburg norm for $\rho_{1,p}$ is equivalent to the Luxemburg norm for the modular obtained by replacing $\rho_p$ with $\eta_{p}$ given by $\eta_p(u)=\int_{\Omega}|u(x)|^{p(x)} dx$. However, it is easy to construct a sequence $(u_j)\subset L^p(\Omega)$ such that $\rho_p(u_j)\rightarrow 0$ but $\eta_p(u_j)\not \rightarrow 0$ as $j\rightarrow \infty$. Indeed, it suffices to take $\Omega=(0,1/2)$, $p(x)= 1/x$ and for each $j\geq 1$, set  $u_j(x)=j^{\frac{2}{j}}\mathbbm{1}_{\left(\frac{1}{j+1},\frac{1}{j}\right)}(x) $.}
		
	\end{remark}
	
	\section{Modular topologies} \label{s4}
	The vast majority of the existing literature on variable exponent Lebesgue spaces on a domain $\Omega \subseteq {\mathbb R}^n$ focuses primarily on the case of bounded variable exponents. The  assumption of boundedness of the exponent $p(x)$ implies that modular convergence and norm convergence are equivalent and that it suffices to concentrate on the Banach space structure of such spaces. There is another feature with far reaching implications if $p(x)$  is bounded away from $1$ and infinity: The condition $1<p_-=\inf\limits_{\Omega}p(x)\leq \sup\limits_{\Omega}p(x)=p_+<\infty$, is equivalent to the Banach space $L^{p(\cdot)}(\Omega)$ being uniformly convex, \cite{Luk}. As a consequence, if the variable exponent is bounded, then geometrically and topologically $L^{p(\cdot)}(\Omega)$ presents close similarities with the classical, constant exponent Lebesgue spaces. On the contrary, if $p(x)$ is unbounded, the Luxemburg norm on $L^{p(\cdot)}(\Omega)$ is not uniformly convex. This lack of smoothness poses serious challenges to the understanding the Banach space structure of the space and stands in the way of the treatment of partial differential equations with non-standard growth in the unbounded case. \\
	Somewhat paradoxically, though, the boundedness assumption obstructs the understanding of the full richness of the modular structure inherent to $L^{p(\cdot)}(\Omega).$\\
	
	As will be evident in the course of this work, the Banach space structure is insufficient to deal with boundary value problems involving differential operators with non-standard growth when the variable exponent is unbounded. It will be shown that the modular topology provides the adequate framework for the treatment of such problems.\\
	Much has been written on modular spaces, see for example \cite{M:1983} and the references therein. An exhaustive treatment of modular functions spaces and their geometry is given in \cite{KK}, Chapters 3 and 4.\\
	The notions presented in what follows concern a general modular vector space $X$; let $\rho:X\rightarrow [0,\infty]$ be a function on a real vector space $X$ having the properties enumerated in (\ref{modularproperties}). Set $X_{\rho}=\{x\in X:\rho(\lambda x)<\infty\,\text{for some}\,\lambda>0\}$. Then \cite{M:1983}, $X_{\rho}$ is a Banach space when furnished with the Luxemburg norm
	\begin{equation*}
		\|x\|_{\rho}=\inf\left\{\lambda>0:\rho\left(\frac{x}{\lambda}\right)\leq 1\right\}.
	\end{equation*}
	Thus, there are two modes of convergence on any modular space: if $(x_j)\subset X_{\rho}$ \cite{KK}
	\begin{enumerate}
		\item [$(i)$] $x_j\overset{\rho}\rightarrow x$ iff $\rho(x_j-x)\rightarrow 0$ as $j\rightarrow \infty$ (modular convergence)\\
		\item[$(ii)$] $x_j\overset{\|\cdot\|_{\rho}}\rightarrow x$ iff $\|x_j-x\|_\rho\rightarrow 0$ as $j\rightarrow \infty$ (norm convergence).
	\end{enumerate}
	Both convergence notions are equivalent if and only if $\rho$ satisfies the $\Delta_2$-condition, i.e., if given any sequence $(x_j)\subset X_{\rho}$ the implication
	\begin{equation}\label{delta2original}
	\rho(x_j)\rightarrow 0\,\text{as}\, j\rightarrow \infty \Rightarrow \rho(2x_j)\rightarrow 0 \,\text{as}\, j\rightarrow \infty
	\end{equation}
	holds. Equivalently, $\rho$ satisfies the $\Delta_2$ condition if and only if there exists a positive constant $K$ such that
	\begin{equation*}\label{delta2}
		\rho(2x)\leq K\rho(x) \,\text{for all}\,x\in X_{\rho}.
	\end{equation*}
	It is well known that the $\Delta_2$-condition holds for the modular space $L^{p}(\Omega)$ if and only if the exponent $p$ is bounded on $\Omega$ \cite{DHHR,KR}.\\
	It is clear then that the modular by itself defines a topology on $X_{\rho}$ (called the modular topology or whose open sets will be characterized shortly) via the modular convergence alluded to above, and it follows from the above discussion that if $p$ is unbounded, the modular topology is different from, in fact, strictly weaker than, the norm topology in any of the modular vector spaces $L^p(\Omega)$ and $W^{1,p}(\Omega)$. Aiming at introducing the solution spaces for the variable exponent $p$-Laplacian, namely $V^{1,p}_0(\Omega)$ and $U^{1,p}_0(\Omega)$, a brief description and some basic properties of the modular topology will be presented next.\\
	
	\begin{definition}\cite{KK}  The open modular ball centered at $x \in X_\rho$ with radius $r > 0$, denoted $B_\rho^o(x,r)$, is defined as
		$$B_\rho^o(x,r) = \{y \in X_\rho;\; \rho(x-y) < r\}.$$
		Similarly, we define the closed modular ball centered at $x \in X_\rho$ with radius $r > 0$, denoted $B_\rho(x,r)$ as
		$$B_\rho(x,r) = \{y \in X_\rho;\; \rho(x-y) \leq r\}.$$
	\end{definition}
	\begin{remark}\label{remarkballlp}{ \normalfont An elementary remark, whose importance will be apparent in the proof of completeness (Theorem \ref{completeness}) is in order here: If $u\in X_{\rho}$, then for $\epsilon>0$ the set 
			$$\left\{v\in X:\rho(v-u)<\epsilon\right\}$$ is entirely contained in $X_{\rho}$. Indeed, for some $k>1$ one has $\rho\left(2^{-k}u\right)<\infty$, so that for $v\in B_\rho(x,\epsilon)$, 
			\begin{align}
				\rho\left(2^{-k-1}v\right)&=\rho\left(2^{-k-1}(v-u)+2^{-k-1}u\right) \\
				& \leq \frac{1}{2}\rho\left(2^{-k}(v-u)\right)
				+\frac{1}{2}\rho\left(2^{-k}u\right)  <\infty. \nonumber
		\end{align}}
	\end{remark}
	The modular balls will play a central role in the definition of the modular topology referred to before.  Modularly closed subsets of a modular space have been introduced in the literature, specifically:
	
	\begin{definition}\cite{KK, M:1983}  The subset $A$ of $X_\rho$ is said to be $\rho$-closed if and only if for any sequence $(x_n)$ in $A$ that $\rho$-converges to $x$, i.e., for which $\lim\limits_{n \to \infty} \rho(x_n-x) = 0$, one has $x \in A$.
	\end{definition}
	
	The terminology is not incidental: The family of all complements of $\rho$-closed subsets of $X_{\rho}$ is in fact a topology on $X_{\rho}$. 
	
	\begin{definition}
		A subset $B\subseteq X_{\rho}$ is said to be $\rho$-open iff and only if $X_{\rho}\setminus B$ is $\rho$-closed.
	\end{definition}
	
	The following result describes the $\rho$-open subsets:
	
	\begin{proposition} Let $A$ be a subset of $X_\rho$.  The following are equivalent
		\begin{enumerate}
			\item A is $\rho$-closed;
			\item for any $x$ in $A^c = X_\rho\setminus A$, there exists $\varepsilon > 0$ such that
			$$B_\rho^o(x,\varepsilon) \subset A^c = X_\rho\setminus A,$$
			i.e., $B_\rho^o(x,\varepsilon) \cap A = \emptyset$.
		\end{enumerate}
	\end{proposition}
	
	\noindent The proof is easy and will be omitted. It is now easy to characterize the $\rho$-open subsets of $X_{\rho}$ via the following proposition, whose proof is elementary.
	
	\begin{proposition}  The subset $B$ of $X_\rho$ is $\rho$-open if and only if for any $x \in B$, there exists $\varepsilon > 0$ such that $B_\rho^o(x,\varepsilon) \subset B$.
	\end{proposition}
	
	It is straightforward to verify that the collection of modularly open subsets of $X_\rho$ is a topology; it will be called the modular topology and denoted by $\tau_\rho$. As is well known, the modularly closed $\rho$-balls are $\rho$-closed provided the modular $\rho$ satisfies the Fatou property, which is the case for the modulars $\rho_p$ and $\rho_{1,p}$. In the sequel, the modular topologies associated to the modulars $\rho_{p}$ and $\rho_{1,p}$  will be denoted by $\tau_p$ and $\tau_{1,p}$, respectively.

	\begin{definition}\cite{KK, M:1983}  Let $A$ be a subset of $X_\rho$.  Define the modular closure of $A$, denoted $\overline{A}^\rho$, as the intersection of all $\rho$-closed subsets of $X_\rho$  containing $A$.
	\end{definition}
	
	\begin{remark} As it is well known (see \cite{KK, M:1983}), we have
		$$A \subset \overline{A} \subset \overline{A}^\rho,$$
		for any subset $A$ in $X_\rho$, where $\overline{A}$ is the closure of $A$ under the topology defined by the Luxemburg norm.  As previously observed, if the modular satisfies the $\Delta_2$-condition, then the two closures coincide. Otherwise, they are generally different. This is the case for $\tau_p$ and $\tau_{1,p}$, when $p$ is unbounded, as is shown in Theorem \ref{Example}.
	\end{remark}
	
	\begin{proposition}\label{basic-properties} The following properties hold:
		\begin{enumerate}
			\item If $U$ is a $\rho$-open subset of $X_\rho$, then $U + x = \{u+x; u \in U\}$ is also $\rho$-open, for any $x \in X_\rho$.  Hence $U+V$ is $\rho$-open provided either $U$ or $V$ is $\rho$-open.
			\item If $U$ is $\rho$-open and $\theta \geq 1$, then $\theta U$ is also $\rho$-open.
			\item For any $x \in \overline{A}^\rho$ and any $\rho$-open subset $U$ such that $x \in U$, then $U\cap A \neq \emptyset$.
			\item If $A$ is convex, then $\overline{A}^\rho$ is convex.
		\end{enumerate}
	\end{proposition}
	\begin{proof}
		\begin{enumerate}
			\item  Let $v \in U+x$, then $v-x \in U$.  Since $U$ is $\rho$-open, there exists $\varepsilon > 0$ such that $B_\rho^o(v-x, \varepsilon) \subset U$.  It is clear that $B_\rho^o(v, \varepsilon) \subset U+x$. As for $U+V$, note that 
			$$U+V = \bigcup_{v \in V} U+v = \bigcup_{u \in U} V+u,$$
			which yields the desired conclusion.
			\item Let $x \in \theta U$.  Then $x/\theta \in U$.  Since $U$ is $\rho$-open, there exists $\varepsilon > 0$ such that $B_\rho^o(x/\theta, \varepsilon) \subset U$.  For any $y \in B_\rho^o(x, \theta \varepsilon)$ it holds
			$$\rho\left(\frac{x}{\theta} - \frac{y}{\theta}\right) \leq \frac{1}{\theta} \rho(x-y) <  \frac{1}{\theta} \ \theta \ \varepsilon = \varepsilon,$$
			i.e., $y/\theta \in B_\rho^o(x/\theta, \varepsilon) \subset U$, which implies $y/\theta \in U$.  Hence $y \in \theta U$, which forces $B_\rho^o(x, \theta \varepsilon) \subset \theta U$.  Therefore, $\theta U$ is $\rho$-open.
			\item Assume that $U\cap A = \emptyset$.  Then  $A \subset U^c = X_\rho\setminus U$.  Since $U$ is $\rho$-open, it follows that $U^c$ is $\rho$-closed.  Hence $\overline{A}^\rho \subset U^c$, i.e., $U \cap \overline{A}^\rho = \emptyset$.  This contradicts the fact that $x$ belongs to $U \cap \overline{A}^\rho$.
			\item Assume that $A$ is convex.  Let $u, v \in \overline{A}^\rho$ and $\alpha \in (0,1)$. If $\alpha u + (1-\alpha)v$ were not in $\overline{A}^\rho$, i.e., if $\alpha u + (1-\alpha)v \in U = X_\rho \setminus \overline{A}^\rho$ it can be easily seen that
			$$u \in U_1 = \frac{1}{\alpha}\ U -\frac{(1-\alpha)}{\alpha}v.$$
			Since $U_1$ is $\rho$-open, it would follow from the previously proved properties  that $A \cap U_1 \neq \emptyset$.  Let $a \in A \cap U_1$. Hence there would exist $u_0 \in U$ such that
			$$a = \frac{1}{\alpha} u_0 - \frac{(1-\alpha)}{\alpha}\ v,$$
			which would imply
			$$v = \frac{1}{(1-\alpha)}\ u_0 - \frac{\alpha}{(1-\alpha)}\ a \in U_2 = \frac{1}{(1-\alpha)}\ U - \frac{\alpha}{(1-\alpha)}\ a.$$
			Hence $A \cap U_2 \neq \emptyset$.  Let $b \in A \cap U_2$.  Since $b \in U_2$, there exists $v_0 \in U$ such that
			$$b = \frac{1}{(1-\alpha)}\ v_0 -\frac{\alpha}{(1-\alpha)}\ a,$$
			which yields $v_0 = \alpha\ a + (1-\alpha)\ b$. It follows from the convexity of $A$ that $\alpha\ a + (1-\alpha)\ b \in A$.  Therefore $v_0 \in A \cap U$ holds.  This contradicts the fact that $U \cap A = \emptyset$.
		\end{enumerate}
	\end{proof}
	
	These properties allow us to prove the following beautiful result:
	
	\begin{proposition}\label{subspace}  Let $A$ be a vector subspace of $X_\rho$.  Then $\overline{A}^\rho$ is a $\rho$-closed vector subspace of $X_\rho$.
	\end{proposition}
	\begin{proof} First, it will be proved that if $u, v \in \overline{A}^\rho$, then $u+v \in \overline{A}^\rho$.  Assume not, i.e., that $u+v \in U = X_\rho\setminus \overline{A}^\rho$.  Hence $u \in (U-v)$, which is $\rho$-open.  This will force $A \cap (U-v) \neq \emptyset$.  Let $a \in A \cap (U-v)$.  There exists $u_0 \in U$ such that $a = u_0 -v$ which yields $v = u_0 - a \in (U-a)$.  Again since $U-a$ is open,
		it follows $A \cap (U-a) \neq \emptyset$.  Let $b \in A \cap (U-a)$.  Hence there exists $v_0 \in U$ such that $b = v_0 - a$ which implies $v_0 = b+a$.  Since $A$ is a subspace, it holds $b+a \in A$.  In other words, the assumption would yield  that $v_0 \in A \cap  U$ which contradicts the fact $A \cap U = \emptyset$.\\
		Next, observe that if $u \in \overline{A}^\rho$ and $\alpha \in \mathbb{R}$, then it must hold that $\alpha\ u \in \overline{A}^\rho$.  Without loss of any generality, it may be assumed that $\alpha \neq 0$.  Let us first prove this conclusion for $\alpha > 0$.  The first part of the proof shows that $k u \in \overline{A}^\rho$, for any $k \in \mathbb{N}$.  So it may be considered that $\alpha$ is not an integer.  In this case, there exists $k \in \mathbb{N}$ such that $k < \alpha < k+1$.  This will imply the existence of $\theta \in (0,1)$ such that $\alpha = \theta \ k + (1-\theta)\ (k+1)$.  Hence
		$$\theta \ u = \theta\ k u + (1-\theta) (k+1) u.$$
		Using the convexity of $A$,  Proposition \ref{basic-properties} yields that $\overline{A}^\rho$ is convex.  Hence $\theta \ u \in \overline{A}^\rho$.  The proof of Proposition \ref{subspace} will be complete by proving that if $u \in \overline{A}^\rho$, then $-u \in \overline{A}^\rho$.  This will follow from identical ideas used before and the fact that if $U$ is $\rho$-open, then $-U$ is also $\rho$-open.  Let's tackle this last statement.  Let $U$ be a $\rho$-open subset of $X_\rho$.  Let $x \in -U$.  Then $-x \in U$ which implies the existence of $\varepsilon > 0$ such that $B_\rho^o(-x, \varepsilon) \subset U$.  Using the properties of the modular it is readily seen that $y \in B_\rho^o(-x, \varepsilon)$ if and only if $-y \in B_\rho^o(x, \varepsilon)$.  Therefore, it holds $B_\rho^o(x, \varepsilon) \subset -U$, which completes the proof that $-U$ is $\rho$-open.
	\end{proof}
	
	\medskip
	
	\begin{definition} \label{defv1p0}
		We denote the $\rho_{1,p}$-closure of $C^{\infty}_0(\Omega)$ in $W^{1,p}(\Omega)$ by $V^{1,p}_0(\Omega)$ and the $\rho_{1,p}$-closure of the subspace of compactly supported functions in $W^{1,p}(\Omega)$ by $U^{1,p}_0(\Omega).$
		In the sequel, $V^{1,p}(\Omega)$ will stand for the set of  all $u \in W^{1,p}(\Omega)$ for which $\rho_{1,p}(u) <\infty$. As usual, the Luxemburg-norm closure of $C_0^{\infty}(\Omega)$ in $W^{1,p}(\Omega)$ will be denoted by $W^{1,p}_0(\Omega)$.
	\end{definition}
	\begin{remark}
		It is obvious that $W^{1,p}_0(\Omega)\subseteq V^{1,p}_0(\Omega)\subseteq U^{1,p}_0(\Omega).$ On account of Theorem 3.11 in \cite{Har}, the strict inclusion $V^{1,p}_0(\Omega)\subsetneq U^{1,p}_0(\Omega)$ might hold even when $p_+<\infty$.
	\end{remark}
	The next result follows directly from Proposition \ref{subspace}.
	
	\begin{theorem}\label{vectorspace}
		$V^{1,p}_0(\Omega)$ and $U^{1,p}_0(\Omega)$ are  $\rho_{1,p}$-closed real vector subspaces of $W^{1,p}(\Omega)$.
	\end{theorem}
	
	\medskip
	\noindent Next, some further results are presented for the specific case of the Lebesgue variable exponent spaces. As discussed above, a subset $A\subset W^{1,p}(\Omega)$ ($A\subset L^{p}(\Omega)$) is $\rho_{1,p}$-closed ($\rho_{p}$-closed) if and only if whenever $(u_j)\subseteq A$ and $\rho_{1,p}(u_j-u)\rightarrow 0$  ($\rho_{p}(u_j-u)\rightarrow 0$) as $j\rightarrow \infty$, it holds that $u\in A$.  If $(u_j)\subset L^p(\Omega)$ ($ W^{1,p}(\Omega)$) and $\lim\limits_{j\rightarrow \infty}\rho_{p}(u_j-u)=0$ ($\lim\limits_{j\rightarrow \infty}\rho_{1,p}(u_j-u)=0$), $(u_j)$ is said to $\rho_p$-converge ($\rho_{1,p}$-converge) to $u$. This will be written as $u_j\overset{\rho_p}{\rightarrow}u$ as $j\rightarrow \infty$ ($u_j\overset{\rho_{1,p}}{\rightarrow}u$ as $j\rightarrow \infty$).\footnote{A different definition of $\rho$ convergence is given in \cite{M:1983}. Our definition is far more stringent, but adequate for the purpose of this work.}\\

	\begin{lemma}\label{aeconvergence}
		If $(u_j)\subseteq L^p(\Omega)$  $\rho_p$-converges to $u\in L^p(\Omega)$, then there exists a subsequence of $(u_j)$ that converges to $u$ a.e. in $\Omega$.
	\end{lemma}
	\begin{proof}
		In \cite[Theorem 12.3, Theorem 12.4]{ML1} the following statement is proved: If $\int_\Omega |f_n-f|\rightarrow 0$ then there exists a subsequence $f_{n_k}$ such that $f_{n_k}\rightarrow f$ a.e.
		Let $\int_\Omega |u_j-u|^p/p\rightarrow 0$. Denoting $v_j=|u_j-u|^p/p$ it is obvious that $\int_\Omega v_j\rightarrow 0$. Applying the previous fact, a subsequence can be chosen such that $v_{j_k}\rightarrow 0$ a.e., from which it follows that $u_{j_k}\rightarrow u$ a.e..
	\end{proof}
	
	As a consequence, one has:
	\begin{corollary}\label{aeconvergence1}
		If $(u_j)\subseteq W^{1,p}(\Omega)$ $\rho_{1,p}$-converges to $u\in W^{1,p}(\Omega)$, then there exists a subsequence $(u_{j_k})$ of $(u_j)$ that converges to $u$ a.e. in $\Omega$ and such that $\left(\nabla u_{j_k}\right)$ converges to $\nabla u$ a.e. in $\Omega.$
	\end{corollary}
	
	\medskip
	The following technical Lemma will be used in the sequel:
	\begin{lemma}\label{stupid}
		The function $f(t)=\displaystyle  \frac{\ln t}{t-1}$ is decreasing on $(1,\infty)$.  Then $f(t)<1$,  which implies $t^{1/(t-1)} < e$, for $t > 1$.
	\end{lemma}
	\begin{proof}
		Let $g(t)=\displaystyle \frac{t-1}{t}-\ln t$. Clearly $g(1)=0$ and
		\begin{align*}
			g'(t)=\frac{1}{t^2}-\frac{1}{t}<0,\ \ \ t\in(1,\infty).
		\end{align*}
		Thus, $g(t)<0$ on $(1,\infty)$ and consequently
		$$f'(t)= \frac{g(t)}{(t-1)^2}<0,\ \ \ t\in(1,\infty).$$
		The fact $f(t)<1$ follows from $\lim\limits_{t\rightarrow 1_+} f(t)=1$.
	\end{proof}
	
	\medskip
	\begin{lemma}\label{pimpliesq}
		Let $(u_j)\subseteq L^p(\Omega)$  $\rho_p$-converge to $u\in L^p(\Omega)$. Assume that $q\leq p$ a.e. and
		\begin{align}
			\int_\Omega \frac{\exp(q(x))}{q(x)} dx = \int_\Omega \frac{e^{q(x)}}{q(x)} dx < \infty.\label{condpq}
		\end{align}
		Then the sequence $(u_j)$ $\rho_q$-converges to $u$.
	\end{lemma}
	\begin{proof}
		Without lost of generality, it can assumed that $u=0$. Hence
		\begin{align*}
			\lim_{j \to \infty} \rho_p(u_j) = \lim_{j \to \infty} \int_\Omega \frac{|u_j(x)|^{p(x)}}{p(x)} dx = 0.
		\end{align*}
		By Lemma \ref{aeconvergence}, there exists a subsequence $u_{j_k}\rightarrow 0$ a.e. Set $\Omega_1= \{x\in\Omega;\ p(x)=q(x)\}$ and $\Omega_2=\Omega\setminus\Omega_1$. Clearly
		\begin{align*}
			&\int_{\Omega_1} \frac{|u_j|^q}{q}=\int_{\Omega_1} \frac{|u_j|^p}{p}\le\int_{\Omega} \frac{|u_j|^p}{p}\rightarrow 0,\ \ \ j\rightarrow\infty.
		\end{align*}
		Define further
		$$A_k=\Big\{x\in\Omega_2; |u_{j_k}|>(p/q)^\frac{1}{p-q}\Big\}\;\; and\;\;\; B_k=\Big\{x\in\Omega_2; |u_{j_k}|\le(p/q)^\frac{1}{p-q}\Big\}.$$
		Let $x\in A_k$. Then it holds $|u_{j_k}|^{p-q}>p/q$, which yields $\displaystyle \frac{|u_{j_k}|^q}{q}\le \frac{|u_{j_k}|^p}{p}$ and consequently
		\begin{align*}
			&\int_{A_k}\frac{|u_{j_k}|^q}{q} \le \int_{A_k}\frac{|u_{j_k}|^p}{p} \le \int_{\Omega}\frac{|u_{j_k}|^p}{p}\rightarrow 0,\ \ \ k\rightarrow\infty.
		\end{align*}
		Consider now $x\in B_{k}$. Denote $\displaystyle \alpha(x)=\frac{p(x)}{q(x)}>1$. Then by Lemma \ref{stupid}
		\begin{align*}
			&|u_{j_k}(x)| \le \left(\frac{p(x)}{q(x)}\right)^{\frac{1}{p(x)-q(x)}} = \alpha(x)^{\frac{1}{q(x)(\alpha(x)-1)}}\le \alpha(x)^{\frac{1}{\alpha(x)-1}}<\e.
		\end{align*}
		By assumption
		\begin{align*}
			&\int_{B_k}\frac{|u_{j_k}|^q}{q} = \int_\Omega\frac{|u_{j_k}\mathbbm{1}_{B_k}|^q}{q} \le \int_{B_k}\frac{\e^{q(x)}}{q(x)} dx \le \int_{\Omega}\frac{\e^{q(x)}}{q(x)} dx < \infty.
		\end{align*}
		Thus one has an integrable majorant and, by virtue of $u_{j_k}\mathbbm{1}_{B_k}\rightarrow 0$ a.e., it follows that
		\begin{align*}
			&\int_{B_k} \frac{|u_{j_k}|^q}{q}\rightarrow 0,\ \ \ k\rightarrow\infty.
		\end{align*}
		Finally,
		$$\lim_{k \to \infty}\ \int_\Omega\frac{|u_{j_k}|^q}{q}= \lim_{k \to \infty}\ \Bigg(\int_{\Omega_1}\frac{|u_{j_k}|^q}{q}+\int_{A_k}\frac{|u_{j_k}|^q}{q} +\int_{B_k}\frac{|u_{j_k}|^q}{q}\Bigg) = 0.$$
		Assume now that the sequence $(u_j)$ does not $\rho_q$-converge to $0$. Then there is $\delta>0$ and a subsequence $(u_{j_k})$ with $\rho_q(u_{j_k})\ge\delta$, for each $k \ge 1$. But the subsequence $(u_{j_k})$ satisfies the assumptions of the lemma; hence there must exist a subsequence $(u_{j_{k_l}})$ such that $\rho_q(u_{j_{k_l}})\rightarrow 0$, which is a contradiction.
	\end{proof}
	
	\begin{corollary}\label{qbouded}
		Let $(u_j)\subseteq L^p(\Omega)$  $\rho_p$-converge to $u\in L^p(\Omega)$. Assume that $q\leq p$ a.e. and $q$ is bounded. Then $(u_j)$ $\rho_q$-converges to $u$.
	\end{corollary}
	
	\medskip
	The above results yield the following:
	
	\begin{theorem}\label{completeness}
		$L^p(\Omega)$ is $\rho_p$-complete, that is, if $(u_j)\subset L^p(\Omega)$ and $\rho_p(u_j-u_k)\rightarrow 0$ as $j,k\rightarrow \infty$, then there exists $u\in L^p(\Omega)$ such that $\rho_p(u_j-u)\rightarrow 0$ as $j\rightarrow \infty$.
	\end{theorem}
	\begin{proof} For any $\varepsilon > 0$, set
		$$\delta(\varepsilon)=\sup\ \left\{\int_\Omega|u(x)|dx;\ \int_\Omega \frac{|u(x)|^{p(x)}}{p(x)}\ dx < \varepsilon\right\}.$$
		Note that the function $\delta$ is increasing.  Hence $\lim\limits_{\varepsilon\rightarrow 0_+}\delta(\varepsilon) = \inf\limits_{\varepsilon > 0}\ \delta(\varepsilon)$.  We claim that
		$\lim\limits_{\varepsilon\rightarrow 0_+}\delta(\varepsilon)=0$. Otherwise, there exist $\eta > 0$ and a sequence $(u_s) \in L^p(\Omega)$ such that
		$$\int_\Omega \frac{|u_s(x)|^{p(x)}}{p(x)} dx < \frac{1}{s} \;\; and \;\;\;  \int_\Omega |u_s(x)|\ dx\ge \eta,$$
		which is in contradiction with Corollary \ref{qbouded}. Thus any $\rho_p$-Cauchy sequence $(u_j)$ must be Cauchy in $L^1(\Omega)$. Then there is $u\in L^1(\Omega)$ which the limit of $(u_j)$ in $L^1(\Omega)$. By Theorems 12.3 and 12.4 in \cite{ML1}, one can select a subsequence $(u_{j_k})$ which converges a.e. to $u$.  Combining the $\rho_p$-Cauchy behavior of $(u_j)$ and the Fatou's lemma it is plain that
		$$\lim_{j \to \infty}\ \rho_p(u-u_j)\leq \lim_{j \to \infty}\ \liminf_{l \to \infty}\ \rho_p(u_{j_l}-u_j) = 0.$$
		Hence, $u$ belongs to an open modular ball centered at $u_j$, and Remark \ref{remarkballlp} yields readily that $u\in L^p(\Omega)$.
		
	\end{proof}
	
	From the above results, the next Theorem  follows immediately.
	\begin{theorem}\label{pqmodularembedding}
		For exponents $q, p \in{\mathcal P}(\Omega)$ with $p\geq q$ a.e. in $\Omega$ satisfying the condition in Lemma \ref{pimpliesq}, the inclusion
		\begin{equation*}
			\tilde{i}_{p,q}: \left(L^p(\Omega), \tau_{p}\right)\longrightarrow \left(L^q(\Omega), \tau_{q}\right)
		\end{equation*}
		is continuous.
	\end{theorem}
	
	\medskip
	Along the same lines, replacing $\rho_p$ with $\rho_{1,p}$, it is easily observed that:
	
	\begin{theorem}\label{pqmodularembedding1}
		Consider the exponent functions $p,q \in{\mathcal P}(\Omega)$,  with $p\geq q$ a.e., with $q$ satisfying the condition in Lemma \ref{pimpliesq}.
		\begin{enumerate}[(i)]
			\item If $(u_j)$ converges to $u$ in $W^{1,p}(\Omega)$, then it must hold $\rho_{1,q}(u_j - u) = 0$, i.e., $(u_j)$ converges to $u $ in $W^{1,q}(\Omega)$.
			\item The inclusion
			\begin{equation*}
				i_{p,q}: \left(W^{1,p}(\Omega), \tau_{1,p}\right)\longrightarrow \left(W^{1,q}(\Omega), \tau_{1,q}\right)
			\end{equation*}
			is continuous.
		\end{enumerate}
	\end{theorem}
	
	\medskip
	A similar result to Theorem \ref{completeness} also holds for $W^{1,p}(\Omega)$.
	
	\begin{theorem}\label{completeness1}
		If $(u_j)\subset W^{1,p}(\Omega)$ and $\rho_{1,p}(u_j-u_k)\rightarrow 0$ as $j,k\rightarrow \infty$, then there exists $u\in W^{1,p}(\Omega)$ such that $\rho_{1,p}(u-u_j)\rightarrow 0$ as $j\rightarrow \infty$. Thus, the Sobolev space $W^{1,p}(\Omega)$ is $\rho_{1,p}$-complete.
	\end{theorem}
	\begin{proof}
		By assumption and Theorem \ref{completeness}, there exist $u\in L^p(\Omega)$ and a vector $\mathbf {v}\in \left(L^p(\Omega)\right)^n$ such that $\rho_{p}(u_j-u)\rightarrow 0$ and $\rho_p\left(\left|\nabla u_j-\mathbf{v}\right|\right)\rightarrow 0$ as $j\rightarrow \infty.$ On account of Theorem \ref{pqmodularembedding1}, the inclusion
		\begin{equation*}
			\left(W^{1,p}(\Omega),\tau_{1,p}\right)\hookrightarrow \left(W^{1,p_-}(\Omega),\tau_{1,p_-}\right)
		\end{equation*}
		is continuous. Therefore the sequence $(u_j)$ is $\rho_{1,p_-}$-Cauchy. Since the norm topology and the $\tau_{1,p_-}$ topology coincide in $W^{1,p_-}(\Omega)$, the sequence $(u_j)$ is Cauchy in the usual Sobolev norm on $W^{1,p_-}(\Omega)$. It follows that $\mathbf{v}=\nabla u$ and therefore $\rho_{1,p}(u_j-u)\rightarrow 0$ as $j\rightarrow \infty.$  As we did in the proof of Theorem \ref{completeness}, we can easily show that $u \in W^{1,p}(\Omega)$.
	\end{proof}
	
	Though $W^{1,p}_0(\Omega)\subseteq V^{1,p}_0(\Omega)$, the next example shows that when $p_+=\infty$ one does not expect equality of both spaces.
	
	\begin{example}\label{Example}  Let $\Omega=(0,1)$ and consider the exponent function $p:(0,1)\rightarrow (0,\infty)$ such that $p(x)=1/x$.  Consider the interval $I_s=(1/(s+1), 1/s]$, for $s \geq 1$. Then, for any $s \geq 1$ and $x \in I_s$, we have $s\leq p(x)<s+1$.  Note that $(0,1]=\bigcup\limits_{s=1}^{\infty} I_s$ holds. Set
		$$w_s(x)=\left(|I_s|^{-1}2^{-s}\right)^{\frac{1}{p(x)}}{\mathbbm{1}}_{I_s}(x) = \left(\frac{s(s+1)}{2^s}\right)^{\frac{1}{p(x)}}{\mathbbm{1}}_{I_s}(x), \; for\;\; s \geq 1,$$
		and $\displaystyle u(x)=\sum\limits_{s =1}^{\infty}\ w_s(x)$.  We claim that the function $\displaystyle v(x)=\int_0^xu(t)dt $ belongs to $V^{1,p}_0(\Omega)$ and not in $W^{1,p}_0(\Omega)$, which will imply $W^{1,p}_0(\Omega) \neq V^{1,p}_0(\Omega)$.\\
		To see this, let us first show that $u$ is integrable.  Indeed,
		\begin{align*}
			\int_0^1u(x)dx&=\sum\limits_{s=1}^{\infty}\int\limits_{I_s}w_s(x)\,dx=
			\sum\limits_{s=1}^{\infty}\int\limits_{I_s}|I_s|^{\frac{-1}{p(x)}}2^{\frac{-s}{p(x)}}dx \\
			&\leq \sum\limits_{s=1}^{\infty}\int\limits_{I_s}|I_s|^{\frac{-1}{s}}2^{\frac{-s}{s+1}}dx\leq \sum\limits_{s=1}^{\infty}|I_s|^{1-\frac{1}{s}}\\
			&= \sum\limits_{s=1}^{\infty}\left(\frac{1}{s(s+1)}\right)^{1-\frac{1}{s}} \leq
			\sum\limits_{s=1}^{\infty}\left(\frac{1}{s^2}\right)^{1-\frac{1}{s}}.
		\end{align*}
		Since
		$$\left(\frac{1}{s^2}\right)^{1-\frac{1}{s}} = \frac{1}{s^{2-\frac{2}{s}}} \leq \frac{1}{s^{1.5}},\; for \;\; s \geq 4,$$
		we conclude that
		$$\int_0^1u(x)dx \leq \sum\limits_{s=1}^{\infty}\left(\frac{1}{s^2}\right)^{1-\frac{1}{s}} < \infty.$$
		Thus, $u$ is integrable, as claimed and $v$ is well defined and bounded. Set
		$$u_k(x) = \sum\limits_{s =1}^k\ w_s(x), \ and \;\; v_k(x) = \int_0^x\ u_k(t)dt,$$
		for $k \geq 1$.  It will be shown next that $\rho_{1,p}\left(v_k-v \right)\rightarrow 0$
		as $k\rightarrow \infty.$ First, notice that
		\begin{align*}
			\rho_p(v^{\prime}_k-v^{\prime}) & = \rho_p(u_k-u) =\rho_p\left(\sum\limits_{s=k+1}^{\infty}w_s\right)\\
			&=\int_0^1 \frac{1}{p(x)} \left(\sum\limits_{s=k+1}^{\infty}w_s(x)\right)^{p(x)}dx =
			\int_0^1\frac{1}{p(x)} \left(\sum\limits_{s=k+1}^{\infty}(w_s(x))^{p(x)}\right) dx\\
			& \leq \int_0^1\sum\limits_{s=k+1}^{\infty}(w_s(x))^{p(x)}dx =\sum\limits_{s=k+1}^{\infty}\int_0^1\left((|I_s|^{-1}2^{-s})^{\frac{1}{p(x)}}\right)^{p(x)}
			\mathbbm{1}_{I_s}(x) dx\\
			& =\sum\limits_{s=k+1}^{\infty}\int_{I_s}\left((|I_s|^{-1}2^{-s})^{\frac{1}{p(x)}}\right)^{p(x)}dx= \sum\limits_{s=k+1}^{\infty}\int_{I_s}|I_s|^{-1}2^{-s}dx \\
			& =\sum\limits_{s=k+1}^{\infty}2^{-s},
		\end{align*}
		which implies $\rho_p(v^{\prime}_k-v^{\prime})\rightarrow 0$ as $k\rightarrow \infty$.  On the other hand, it follows from the integrability of $u(x)$ that
		\begin{align*}
			\lim\limits_{k \to \infty}\ \int_0^x \Big(u(t)-u_k(t)\Big)dt &\leq \lim\limits_{k \to \infty}\ \int_0^1 \Big(u(t)-u_k(t)\Big) dt\\
			&\leq \lim\limits_{k \to \infty}\  \int_0^1 u(t)\ \mathbbm{1}_{[0, 1/k+1]}(t) dt = 0.
		\end{align*}
		Thus, for sufficiently large $k$, $\displaystyle \int_0^x \Big(u(t)-u_k(t)\Big) dt\leq 1$, for all $x\in (0,1)$, and
		\begin{equation*}
			\left(\int_0^x (u(t)-u_k(t))dt\right)^{p(x)}\leq 1.
		\end{equation*}
		It follows easily from Lebesgue's dominated convergence that
		\begin{equation*}
			\lim\limits_{k \to \infty}\ \rho_p(v_k-v) = 0.
		\end{equation*}
		Thus, $v$ is the modular limit of a sequence of functions in $W^{1,p}(\Omega)$ with compact support in $(0,1]$. It will be shown next that $v$ is the modular limit of a sequence in $C_0^{\infty}(\Omega)$. To this end, let $\varepsilon \in (0, 1/2)$ and select $\delta>0$ and $k\in {\mathbb N}$, such that $1/(k+1) < \delta$, satisfying the following conditions:
		\begin{enumerate}[(i)]
			\item $\displaystyle \int_0^{\delta}\left(|v(t)|^{p(t)}+|v^{\prime}(t)|^{p(t)}\right)dt <\varepsilon$
			\item $\displaystyle \sum\limits_{s=k+1}^{\infty}|I_s|^{1-\frac{1}{s}}<\varepsilon$;
			\item $\displaystyle \int_0^1\left(|v(x)-v_k(x)|^{p(x)} + |v^{\prime}(x)-v^{\prime}_k(x)|^{p(x)}\right)dx < \varepsilon.$
		\end{enumerate}
		Let $\lambda=2^{1-1/\delta}\varepsilon<1.$ On account of Theorem 2 in \cite{ER}, there exists a function $\varphi_k$ in $C^{\infty}_0(\Omega)$ such that
		\begin{enumerate}
			\item $(\text{supp}\ \varphi_k) \cap (0,\delta)=\emptyset,$
			\item $\|v-\varphi_k\|_{1,p}<\lambda<1.$
		\end{enumerate}
		Condition $(2)$ above yields
		\begin{equation*}
			\int_0^1 \left(\left|\frac{v_k(x)-\varphi_k(x)}{\lambda}\right |^{p(x)}+\left|\frac{v_k^{\prime}(x)-\varphi^{\prime}_k(x)}{\lambda}\right|^{p(x)}\right)dx<1.
		\end{equation*}
		Since $\lambda < 1$, we get
		$$\int_{\delta}^1\frac{|v_k(x)-\varphi_k(x)|^{p(x)}}{\lambda} dx  \leq \int_{\delta}^1 \left|\frac{v_k(x)-\varphi_k(x)}{\lambda}\right|^{p(x)} dx$$
		and
		$$\int_{\delta}^1 \frac{|v_k^{\prime}(x)-\varphi^{\prime}_k(x)|^{p(x)}}{\lambda} dx \leq
		\int_{\delta}^1 \left|\frac{v_k^{\prime}(x)-\varphi^{\prime}_k(x)}{\lambda}\right |^{p(x)} dx.$$
		By virtue of $(i)$ and $(1)$ it follows that
		$$\int_{\delta}^1 \left(\frac{|v_k(x)-\varphi_k(x)|^{p(x)}}{\lambda} +\frac{|v_k^{\prime}(x)-\varphi^{\prime}_k(x)|^{p(x)}}{\lambda} \right)dx < 1.$$
		Thus,
		$$\int_{\delta}^1|v_k(x)-\varphi_k(x)|^{p(x)} +|v_k^{\prime}(x)-\varphi^{\prime}_k(x)|^{p(x)}dx < \lambda  = 2^{1-1/\delta}\varepsilon. $$
		The final step is to estimate
		$$\rho_{1,p}(v-\varphi_k) \leq \int_0^1|v-\varphi_k|^{p(x)}dx + \int_0^1|v^{\prime}-\varphi^{\prime}_k|^{p(x)}dx.$$
		Set
		$$B_1 = \int_0^1|v-\varphi_k|^{p(x)}dx \;\;\; and \;\;\; B_2 = \int_0^1|v^{\prime}-\varphi^{\prime}_k|^{p(x)}dx.$$
		As for $B_1$, we have
		\begin{align*}
			B_1&=\int_0^1|v-\varphi_k|^{p(x)}dx=\int_0^{\delta}|v-\varphi_k|^{p(x)}dx+\int_{\delta}^1|v-\varphi_k|^{p(x)}dx\\
			&=\int_0^{\delta}|v|^{p(x)}dx+\int_{\delta}^1|v-v_k+v_k-\varphi_k|^{p(x)}dx\\
			&=\int_0^{\delta}|v|^{p(x)}dx + \int_{\delta}^1\left|2 \frac{v-v_k+v_k-\varphi_k}{2}\right|^{p(x)}dx\\
			&\leq \varepsilon+\int_{\delta}^12^{p(x)-1}|v-v_k|^{p(x)}dx+
			2^{p(\delta)-1}\int_{\delta}^1|v_k-\varphi_k|^{p(x)}dx\\
			& \leq 2\varepsilon+ \int_{\delta}^1 2^{p(x)-1}|v-v_k|^{p(x)}dx.
		\end{align*}
		On the other hand, since $1/(s+1) < 1/p(x) \leq 1/s$, for any $x \in I_s$, we get
		\begin{align*}
			v(x)-v_k(x)&=\int_0^x(u(t)-u_k(t))dt\leq \int_0^1 \Big(\sum\limits_{s=k+1}^{\infty}w_s(t)\Big) dt\\
			&= \int_0^1\sum\limits_{s=k+1}^{\infty}\left(|I_s|^{-1}2^{-s}\right)^{\frac{1}{p(t)}} \mathbbm{1}_{I_s}(t) dt = \sum\limits_{s=k+1}^{\infty}\int_{I_s} \left(|I_s|^{-1}2^{-s}\right)^{\frac{1}{p(t)}} dt \\
			& \leq \sum\limits_{s=k+1}^{\infty}\int_{I_s}|I_s|^{\frac{-1}{s}}2^{\frac{-s}{s+1}}dt\\
			& \leq \sum\limits_{s=k+1}^{\infty}|I_s|^{1-\frac{1}{s}}\leq \varepsilon,
		\end{align*}
		for any $x \in [0,1]$.  It follows
		\begin{align*}
			\int_0^12^{p(x)-1}(v(x)-v_k(x))^{p(x)}dx&\leq \int_0^12^{p(x)-1}\varepsilon^{p(x)}dx =\frac{1}{2}\int_0^1(2\varepsilon)^{p(x)}dx \\ \nonumber &\leq \frac{1}{2}\int_0^12\varepsilon \, dx=\varepsilon,
		\end{align*}
		and hence $B_1\leq 3\varepsilon$.  As to $B_2$, one has
		\begin{align*}
			B_2&=\int_0^1(v^{\prime}-\varphi_k^{\prime})^{p(x)}dx=\int_0^{\delta}(v^{\prime}- \varphi_k^{\prime})^{p(x)}dx +\int_{\delta}^1(v^{\prime}-\varphi_k^{\prime})^{p(x)}dx\\ &=\int_0^{\delta}(v^{\prime})^{p(x)}dx + \int_{\delta}^1 2^{p(x)}\left(\frac{v^{\prime}-v_k^{\prime}+v_k^{\prime} -\varphi_k^{\prime}}{2}\right)^{p(x)}dx \\
			& \leq \varepsilon+2^{p(\delta)}\left(\frac{1}{2}\int_{\delta}^1(v^{\prime} -v_k^{\prime})^{p(x)}dx + \frac{1}{2}\int_{\delta}^1(v_k^{\prime}-\varphi_k^{\prime})^{p(x)} dx \right)\\
			& \leq \varepsilon +2^{p(\delta)-1}\int_{\delta}^1(v^{\prime}-v_k^{\prime})^{p(x)}dx+\varepsilon.
		\end{align*}
		If $x>\delta$, then $x> 1/(k+1)$, whence
		$$v^{\prime}(x)-v_k^{\prime}(x)=u(x)-u_k(x)=\sum\limits_{s=k+1}^{\infty}w_s(x)=0,$$
		and consequently $B_2\leq 2\varepsilon.$ It is clear from the above estimates that $\varrho_{1,p}(v-\varphi_k)<5\varepsilon$.  Therefore, we have $v\in V^{1,p}_0(\Omega)$.  In order to show that $v\notin W^{1,p}_0(\Omega)$, consider $\varphi$ in $C^{\infty}_0(\Omega)$ with $\text{supp} (\varphi) \cap (0,1/(k+1))=\emptyset.$ Then
		\begin{align*}
			\rho_{1,p}\Big(2(v-\varphi)\Big) & \geq \int_0^1 \frac{1}{p(x)} \Big(2(v^{\prime}-\varphi^{\prime})\Big)^{p(x)}dx \geq \int_0^{1/(k+1)} \frac{1}{p(x)} 2^{p(x)}|v^{\prime}|^{p(x)} dx \\
			&= \sum\limits_{s=k+1}^{\infty} \int_{I_s}\frac{1}{p(x)} 2^{p(x)}|u|^{p(x)}dx =\sum\limits_{s=k+1}^{\infty}\int_{I_s} \frac{1}{p(x)} 2^{p(x)}|w_s|^{p(x)}dx \\
			&= \sum\limits_{s=k+1}^{\infty}\int_{I_s} \frac{1}{p(x)} 2^{p(x)}|I_s|^{-1}2^{-s} dx \geq
			\sum\limits_{s=k+1}^{\infty}\int_{I_s} \frac{1}{s+1} 2^s|I_s|^{-1}2^{-s} dx\\  &=\sum\limits_{s=k+1}^{\infty}\frac{1}{s+1} = \infty.
		\end{align*}
		By definition of the Luxemburg norm, it must hold
		$$\|v-\varphi\|_{1,p}\geq \frac{1}{2},$$
		which proves our claim.
	\end{example}
	
	\medskip\medskip
	The next lemma is evidence of the fact that $V^{1,p}_0(\Omega) \subset U^{1,p}_0(\Omega)$ is not an artificial construction, but it seems to be the adequate substitute of the classical Sobolev space $W^{1,p}_0(\Omega)$.
	
	\begin{lemma}
		Consider $p\in {\mathcal P}(\Omega)$.
		\begin{itemize}
			\item [(i)] If $u\in  U^{1,p}_0(\Omega)\cap C(\overline{\Omega})$, then $u$ vanishes on $\partial\Omega$. \\
			\item[(ii)] If $u\in V^{1,p}(\Omega)\cap C(\overline{\Omega})$ and $u$ vanishes on $\partial \Omega$, then $u\in U^{1,p}_0(\Omega).$
		\end{itemize}
	\end{lemma}
	\begin{proof}
		$(i)$ Observe that for any $\delta>0$, the set
		$$U=\{v\in W^{1,p}(\Omega):\|u-v\|_{L^1(\Omega)}+\||\nabla(u-v)|\|_{L^1(\Omega)}<\delta\}$$
		is a $\rho_{1,p}$-open subset of $W^{1,1}(\Omega)$ containing $u$. To see this, consider a sequence $(v_j)$ in $W^{1,p}(\Omega)\setminus U$ which $\rho_{1,p}$ converges to $v \in W^{1,p}(\Omega)$.  Let us show that $v \in W^{1,p}(\Omega)\setminus U$.  Indeed, virtue of Theorem \ref{pqmodularembedding1}, $(v_j)$ norm converges to $v$ in $W^{1,1}(\Omega)$.  Since
		\begin{align*}
			&\int_{\Omega}|v_j-u|\ dx + \int_{\Omega}|\nabla (v_j-u)|\ dx \ge \delta,\\
			&\lim_{j \to \infty} \Big(\int_{\Omega}|v_j-v|\ dx + \int_{\Omega}|\nabla (v_j-v)|\ dx\Big) = 0,\\
			&\int_{\Omega}|v_j-u|dx - \int_{\Omega}|v_j-x|dx \le \int_{\Omega}|v-u|dx, \\
			&\int_{\Omega}|\nabla (v_j-u)|dx - \int_{\Omega}|\nabla (v_j-v)|dx \leq \int_{\Omega}|\nabla (v-u)|dx,
		\end{align*}
		we get
		$$\int_{\Omega}|v-u|dx +\int_{\Omega}|\nabla (v-u)|dx \ge \delta,$$
		i.e., $v \in W^{1,p}(\Omega)\setminus U$.  It therefore follows that for $u\in U^{1,p}_0(\Omega)$ there is a sequence $(u_j)\subset W^{1,p}(\Omega)$ (each $u_j$ compactly supported) that converges to $u$ in $W^{1,1}(\Omega)$.\\
		Let $P \in \partial \Omega$ be a point on the boundary. After rotations and translations, the boundary of $\Omega\subset {\mathbb R}^n$ can be locally described in a neighborhood $N(P)$ of $P$ as the graph of a $C^2$ function  $\psi: {\mathbb R}^{n-1}\rightarrow {\mathbb R}$, namely as $(x,\psi(x))$, for $|x|\leq 1$, $x\in {\mathbb R}^{n-1}$. In this local coordinate system, any point in $\Omega \cap N(P)$ can be described as $(x,\psi(x)+t)$, $t>0$. In this notation, for each $j$ and any $x\in {\mathbb R}^{n-1}$ such that $|x|\leq 1$, it holds
		\begin{equation*}
			|u_j(x,\psi(x)+t)|\leq \int_0^{t}\left|\frac{\partial u_j}{\partial x_n}(x,\psi(x)+s)\right|ds.
		\end{equation*}
		Integrating the preceding inequality on $(0,\delta)$, $\delta>0$, it is clear that
		\begin{equation*}
			\int_0^{\delta}|u_j(x,\psi(x)+t)|dt\leq \int_0^{\delta}\int_0^{t}\left|\frac{\partial u_j}{\partial x_n}(x,\psi(x)+s) \right| ds\ dt\leq \delta \int_0^{\delta}\left|\frac{\partial u_j}{\partial x_n}(x,\psi(x)+s)\right|ds.
		\end{equation*}
		Hence
		\begin{equation*}\label{intuj}
			\int\limits_{|x|\leq 1}\int_0^{\delta}|u_j(x,\psi(x)+t)|dt\ dx \leq \delta \int\limits_{|x|\leq 1} \int_0^{\delta}\left|\frac{\partial u_j}{\partial x_n}(x,\psi(x)+s)\right|ds\ dx.
		\end{equation*}
		Because $u_j\rightarrow u$ in $W^{1,1}(\Omega)$, letting $j\rightarrow \infty$ on both sides of the above inequality yields
		\begin{equation*}
			\int\limits_{|x|\leq 1}\int_0^{\delta}|u(x,\psi(x)+t)|dt\ dx \leq \delta \int\limits_{|x|\leq 1} \int_0^{\delta}\left|\frac{\partial u}{\partial x_n}(x,\psi(x)+s)\right|ds\ dx,
		\end{equation*}
		which implies
		\begin{equation*}
			\int\limits_{|x|\leq 1}\frac{1}{\delta}\int_0^{\delta}|u(x,\psi(x)+t)|dt\ dx \leq  \int\limits_{|x|\leq 1} \int_0^{1}\left|\frac{\partial u}{\partial x_n}(x,\psi(x)+s)\right|{\mathbbm{1}}_{(0,\delta)}(s)ds\ dx.
		\end{equation*}
		Since $u\in C(\overline{\Omega})$, it is clear, if we let $\delta\rightarrow 0$, that $u$ must vanish on $\partial \Omega$, as claimed.\\
		
		\noindent $(ii)$ Assume now that $u\in C(\overline{\Omega})\cap V^{1,p}(\Omega)$ and $u(x)=0$ for each $x\in \partial \Omega$. Define $u_+(x):=\max \{u(x),0\}$ and $u_-(x):=\min \{u(x),0\}$.  For any $\varepsilon>0$, set
		$$u_{+,\varepsilon}(x):=\max\{u_+(x)-\varepsilon,0\}\;\; and \;\; u_{-,\varepsilon}(x):= \min\{u_-(x)+\varepsilon,0\}.$$
		Then write $u_{\varepsilon}(x):= u_{+,\varepsilon}(x) + u_{-,\varepsilon}(x)$.  It follows easily that $u_\varepsilon \in U^{1,p}_0(\Omega)$.
		Obviously $\rho_{1,p}(u_{\varepsilon}) \le \rho_{1,p}(u),$
		$|\supp u \setminus \supp u_{\varepsilon} | \searrow 0$ as $\varepsilon \searrow 0$, and $\nabla u =\nabla u_{\varepsilon}$ on $\supp u_\varepsilon$. Using this and the fact that $u, u_{\varepsilon} \in C(\overline{\Omega})$, it is readily concluded that $\rho_{1,p}(u-u_{\varepsilon})\to 0$ as $\varepsilon \to 0$, from which follows $u\in U^{1,p}_0(\Omega)$.
	\end{proof}
	
	\medskip
	\begin{lemma}
		If $v\in V^{1,p}_0(\Omega)$ (resp. $v \in U^{1,p}_0(\Omega)$) and $\phi\in C_0^{\infty}(\Omega)$, then $v\phi\in V^{1,p}_0(\Omega)$ (resp. $v\phi \in U^{1,p}_0(\Omega)$).
	\end{lemma}
	\begin{proof}
		Let $\mathcal{O}$ be $\rho_{1,p}$-open in $W^{1,p}(\Omega)$ with $v\phi \in \mathcal{O}$. Let $A=\{u\in W^{1,p}(\Omega):u\phi \in \mathcal{O}\}.$ Clearly $v\in A$. We claim that $A$ is $\rho_{1,p}$-open. To this end, let $(a_j)\subseteq W^{1,p}(\Omega)\setminus A$ and assume that $a_j\overset{\rho_{1,p}}\longrightarrow a.$ If $a\phi\in \mathcal{O}$, then for some $\delta>0$ the modular ball $B_{\rho_{1,p}}(a\phi, \delta)$ must be contained in $\mathcal{O}$, since $\mathcal{O}$ is $\rho_{1,p}$-open. Then, since
		\begin{equation*}
			\int_{\Omega} \frac{|a\phi-a_j\phi|^{p(x)}}{p(x)}\ dx + \int_{\Omega} \frac{|\nabla (a\phi-a_j\phi)|^{p(x)}}{p(x)} \ dx \rightarrow 0, \;\;\; as\;\; j\rightarrow \infty,
		\end{equation*}
		it would follow that $a_i\phi\in B^{\rho_{1,p}}_{\delta}(a\phi)\subset \mathcal{O}$ for some large enough $i$. In other words, the assumption would imply that $a_i\in A$, which contradicts the choice of $(a_i)$. In all, $A$ is $\rho_{1,p}$-open, as claimed, and contains $v$.\\
		Next, assuming $v\in V^{1,p}_0(\Omega)$, there must exist $\psi\in C^{\infty}_0(\Omega)\cap A$. By definition, one then has $\psi \phi\in \mathcal{O}$. Clearly $\psi \phi\in C^{\infty}_0(\Omega)$. Thus, any $\rho_{1,p}$-open set $\mathcal{O}$ containing $v\phi$ intersects $C^{\infty}_0(\Omega)$, which by definition means $v\phi \in V^{1,p}_0(\Omega).$ On the other hand, if $v\in U^{1,p}_0(\Omega)$, one can find  a compactly supported function $\psi \in W^{1,p}(\Omega)\cap A$. Then  $\psi\phi\in \mathcal{O}$ and since $\psi\phi$ is compactly supported and belongs to $W^{1,p}(\Omega)$, it follows that $v\psi\in U^{1,p}_0(\Omega)$.
	\end{proof}
	
	\section{Modular inequalities}\label{s5}
	
	With specific uniform-convexity like modular inequalities on the Sobolev space $W^{1,p}(\Omega)$ in sight, this section starts with the following Lemma.
	
	\begin{lemma} \label{vector-p-inequalities}  Let ${\mathbf u}$, ${\mathbf v}$ be vectors in a Hilbert space $({\mathbb H}, \|\cdot\|)$.  If $1\leq p\leq 2$ it holds
		\begin{equation}\label{1<p<2-vector}
			\left\|\frac{\mathbf{u} + \mathbf{v}}{2}\right\|^p + \frac{p(p-1)}{2^{p+1}}\frac{\|\mathbf{u}-\mathbf{v}\|^2}{(\|\mathbf{u}\|+\|\mathbf{v}\|)^{2-p}} \leq \frac{1}{2}(\|\mathbf{u}\|^p+\|\mathbf{v}\|^p),
		\end{equation}
		provided $\|\mathbf u\| + \|\mathbf v\| \neq 0$.  In addition, if $p\geq 2$ it holds
		\begin{equation}\label{p>2-vector}
			\left\|\frac{\mathbf{u} + \mathbf{v}}{2}\right\|^p + \left\|\frac{\mathbf{u} - \mathbf{v}}{2}\right\|^p \leq \frac{1}{2}(\|\mathbf{u}\|^p+\|\mathbf{v}\|^p).
		\end{equation}
	\end{lemma}
	\begin{proof}
		See \cite{KaMe}.
	\end{proof}
	\begin{corollary} \label{gradient-p-inequalities}  Let $u$, $v$ be functions in $W^{1,p}(\Omega)$ and $|\cdot|$ stand for the Euclidean norm in ${\mathbb R}^n$.  If $1\leq p\leq 2$ it holds
		\begin{equation}\label{1<p<2-vector}
			\frac{1}{p}\left|\frac{\nabla u + \nabla v}{2}\right|^p + \frac{(p-1)}{2^{p+1}}\frac{|\nabla u-\nabla v|^2}{(|\nabla u|+|\nabla v|)^{2-p}} \leq \frac{1}{2p}(|\nabla u|^p+|\nabla v|^p),
		\end{equation}
		provided $|\nabla  u| + |\nabla v| \neq 0$.  In addition, if $p\geq 2$ it holds
		\begin{equation}\label{p>2-vector}
			\frac{1}{p}\left|\frac{\nabla u + \nabla v}{2}\right|^p + \frac{1}{p}\left|\frac{\nabla u - \nabla v}{2}\right|^p \leq \frac{1}{2p}(|\nabla u|^p+|\nabla v|^p).
		\end{equation}
	\end{corollary}
	The following definition is in order, since it will have an essential place in the sequel (see \cite{DHHR} for further details)
	\begin{definition}\label{defuc}
		A convex modular $\rho$ on a vector space $V$ is said to be uniformly convex (in short $(UC^*)$) if for every $\varepsilon>0$ there exists $\delta=\delta(\varepsilon)>0$ such that  for every $u\in V$ and $v\in V$:
		$$\rho\left(\frac{u-v}{2}\right) \geq \varepsilon\ \frac{\rho(u)+\rho(v)}{2} \;\; \text{implies}\;\;\;
		\rho\left(\frac{u+v}{2}\right) \leq (1-\delta)\ \frac{\rho(u)+\rho(v)}{2}.$$
	\end{definition}
	
	\begin{theorem}\label{N=2}
		\label{UC-sobolev2} Let $\Omega\subseteq {\mathbb R}^n$ be a domain, $p\in {\mathcal P}(\Omega)$ finite a.e. and (here for $u\in W^{1,p}(\Omega), u_j=\frac{\partial u}{\partial x_j}$ and $|\cdot|$ stands for the Euclidean norm in ${\mathbb R}^n$).  Consider the functional $\varrho :W^{1,p}(\Omega)\rightarrow [0,\infty]$ defined by
		$$\varrho(u) := \rho_{p}(|\nabla u|)= \int_{\Omega}\frac{1}{p(x)} \left(\sum_1^nu_j^2\right)^{\frac{p(x)}{2}}\ dx = \int_{\Omega}\frac{1}{p(x)}|\nabla u|^{p(x)}\ dx.$$
		Then $\varrho$ is a convex pseudomodular (i.e., it has all the properties exhibited in Definition \ref{def-rho}, except $(1)$) on $W^{1,p}(\Omega)$ and is $(UC^*)$ provided $p_->1$.
	\end{theorem}
	\begin{proof}
		It is obvious that $\varrho$ is a convex pseudomodular on $W^{1,p}(\Omega)$. It remains to prove the uniform convexity assertion. Without of generality, we assume $\max\{\varrho (u),\varrho(v)\}$ is finite.  In the course of the proof, it will be understood that $|(x_1,x_2,...,x_n)|=\left(\sum\limits_1^nx_j^2\right)^{\frac{1}{2}}$ and for a subset $A\subseteq \Omega$, we set $\varrho_{A}( u)=\int_A p^{-1}|{\nabla u}|^p\ dx$.  Let $\varepsilon \in (0,1)$.  Assume
		$$\varrho \left(\frac{ u- v}{2}\right)\geq \varepsilon \ \frac{\varrho( u)+\varrho( v)}{2}.$$
		Set $\Omega_1=\{x\in \Omega: p(x)\geq 2\}$; then necessarily, either
		\begin{equation}\label{omega1}
			\int_{\Omega_1}\frac{1}{p}\left|\frac{\nabla u-\nabla  v}{2}\right|^{p}\ dx \geq \frac{\varepsilon}{2}\left(\frac{\varrho(  u)+\varrho( v)}{2}\right)
		\end{equation}
		or
		\begin{equation}\label{minusomega1}
			\int_{\Omega\setminus \Omega_1}\frac{1}{p}\left|\frac{\nabla u-\nabla v}{2}\right|^{p}\ dx \geq \frac{\varepsilon}{2}\left(\frac{\varrho( u)+\varrho( v)}{2}\right).
		\end{equation}
		In the first case, by virtue of Corollary \ref{gradient-p-inequalities}, it is readily concluded that
		\begin{align*}
			\int_{\Omega_1}\frac{1}{p}\left|\frac{\nabla u+\nabla  v}{2}\right|^{p}\ dx&\leq \frac{1}{2}\left( \int_{\Omega_1}\frac{1}{p}|\nabla u|^p\,dx+\int\limits_{\Omega_1}\frac{1}{p}|\nabla v|^p\,dx\right)- \frac{\varepsilon }{2}\left(\frac{\varrho(u)+\varrho( v)}{2}\right)\\
			&\leq \frac{1}{2}\left( \int\limits_{\Omega_1}\frac{1}{p}|\nabla  u|^p\,dx+\int_{\Omega_1}\frac{1}{p}|\nabla v|^p\,dx\right)- \frac{\varepsilon }{2}\left(\frac{\varrho_{\Omega_1}( u)+\varrho_{\Omega_1}( v)}{2}\right).
		\end{align*}
		In all,
		\begin{align*}\nonumber
			\varrho\left(\frac{ u+ v}{2}\right) &=\int_{\Omega_1}\frac{1}{p}\left|\frac{\nabla u+\nabla v}{2}\right|^{p}\,dx+\int_{\Omega \setminus\Omega_1}\frac{1}{p}\left|\frac{\nabla u+\nabla v}{2}\right|^{p}\,dx\\ \nonumber &\leq
			\frac{1}{2}\left( \int_{\Omega_1}\frac{1}{p}|\nabla u|^p\,dx+\int_{\Omega_1}\frac{1}{p}|\nabla v|^p\,dx\right)- \frac{\varepsilon }{2}\left(\frac{\varrho( u)+\varrho(  v)}{2}\right)
			\\ \nonumber &+\int_{\Omega \setminus\Omega_1}\frac{1}{p}\left|\frac{\nabla u+\nabla v}{2}\right|^{p}\,dx\\ \nonumber &\leq
			\frac{1}{2}\left( \int_{\Omega_1}\frac{1}{p}|\nabla u|^p\,dx+\int_{\Omega_1}\frac{1}{p}|\nabla v|^p\,dx\right)- \frac{\varepsilon }{2}\left(\frac{\varrho( u)+\varrho( v)}{2}\right)
			\\ \nonumber &+\frac{1}{2}\left(\int_{\Omega \setminus\Omega_1}\frac{1}{p}|\nabla u|^p\,dx+\int_{\Omega\setminus \Omega_1}\frac{1}{p}|\nabla v|^{p}\,dx\right)\\ \nonumber &=
			\frac{1}{2}(\varrho( u)+\varrho( v))-\frac{\varepsilon }{2}\left(\frac{\varrho( u)+\varrho( v)}{2}\right)\\ \nonumber &=\left(1-\varepsilon/2\right)\left(\frac{\varrho( u)+\varrho( v)}{2}\right).
		\end{align*}
		The last statement settles the issue in case (\ref{omega1}) holds.  If instead (\ref{minusomega1}) holds, define
		$$\Omega_2=\left\{x\in \Omega\setminus \Omega_1:|\nabla u-\nabla v|<\frac{\varepsilon}{4}\left(|\nabla u|^2+|\nabla  v|^2\right)^{\frac{1}{2}}\right\}.$$
		It follows
		\begin{align} \nonumber
			\int_{\Omega_2}\frac{1}{p}\left|\frac{\nabla u-\nabla v}{2}\right|^{p}\,dx&\leq \int_{\Omega_2}\frac{1}{p}\left|\frac{\varepsilon}{4}\frac{1}{2}(|\nabla u|^2+|\nabla v|^2)^{\frac{1}{2}}\right|^p\,dx \\ \nonumber &\leq \frac{\varepsilon}{4}\int_{\Omega_2}
			\frac{1}{p}\left|\frac{1}{2}(|\nabla u|+|\nabla v|)\right|^p\,dx\\ \nonumber &\leq \frac{\varepsilon}{8}\left(\varrho( u)+\varrho( v)\right).
		\end{align}
		On account of inequality (\ref{minusomega1}), it is readily obtained that
		\begin{align*}\label{geqto}
			\int_{\Omega\setminus (\Omega_1\cup \Omega_2)}\frac{1}{p}\left|\frac{\nabla u-\nabla v}{2}\right|^{p}\,dx& = \int_{\Omega\setminus \Omega_1}\frac{1}{p}\left|\frac{\nabla u-\nabla v}{2}\right|^{p}\,dx - \int_{\Omega_2}\frac{1}{p}\left|\frac{\nabla u-\nabla v}{2}\right|^{p}\,dx\\ \nonumber &\geq \frac{\varepsilon }{2}\left(\frac{\varrho( u)+\varrho( v)}{2}\right)-\frac{\varepsilon}{4}\left(\frac{\varrho( u)+\varrho( v)}{2}\right)\\ \nonumber&=\frac{\varepsilon}{4}\left(\frac{\varrho( u)+\varrho(v)}{2}\right)
		\end{align*}
		By definition, one has, on $\Omega\setminus \left(\Omega_1\cup \Omega_2\right)$:
		\begin{align}\nonumber
			\left|\frac{\nabla u+\nabla v}{2}\right|^p
			+(p_--1)\frac{\varepsilon}{8}\left|\frac{\nabla u-\nabla v}{2}\right|^p &\leq \left|\frac{\nabla u+\nabla v}{2}\right|^p+
			\frac{p(p-1)}{2}\left(\frac{\varepsilon}{4}\right)^{2-p}\left|\frac{\nabla u-\nabla v}{2}\right|^p
			\\ \nonumber &\leq \left|\frac{{\nabla u}+{\nabla v}}{2}\right|^p
			\\ \nonumber &+
			\frac{p(p-1)}{2}\left(\frac{|{\nabla u}-{\nabla v}|}{(|{\nabla u}|^2+|{\nabla v}|^2)^{\frac{1}{2}}}\right)^{2-p}\left|\frac{{\nabla u}-{\nabla v}}{2}\right|^p \\ \nonumber &
			\leq \frac{1}{2}\left(|{\nabla u}|^p+|{\nabla v}|^p \right),
		\end{align}
		on account of the first part of Corollary \ref{gradient-p-inequalities}. Taking into consideration (\ref{minusomega1}), it follows that
		\begin{align}\nonumber
			\int_{\Omega\setminus \Omega_1\cup \Omega_2}\frac{1}{p}\left|\frac{{\nabla u}+{\nabla v}}{2}\right|^p\,dx&\leq \frac{1}{2}\left(\int_{\Omega\setminus \Omega_1\cup \Omega_2}\frac{1}{p}|{\nabla u}|^p\,dx+\int_{\Omega\setminus \Omega_1\cup \Omega_2}\frac{1}{p}|{\nabla v}|^p \right)\\ \nonumber &-(p_--1)\frac{\varepsilon^2 }{32}\left(\frac{\varrho({ u})+\varrho({ v})}{2}\right).
		\end{align}
		Finally,
		\begin{align*}
			\varrho\left(\frac{{ u}+{ v}}{2}\right) &= \int_{\Omega_1\cup\Omega_2} \frac{1}{p}\left|\frac{{\nabla u}+{\nabla v}}{2}\right|^p\,dx + \int_{\Omega\setminus(\Omega_1\cup\Omega_2)} \frac{1}{p}\left|\frac{{\nabla u}+{\nabla v}}{2}\right|^p\,dx\\
			&\leq \frac{1}{2}\ \int_{\Omega_1\cup \Omega_2}\frac{1}{p}|{\nabla u}|^p\,dx + \frac{1}{2} \int_{\Omega_1\cup \Omega_2}\frac{1}{p}|{\nabla v}|^p\,dx + \frac{1}{2} \int_{\Omega\setminus \Omega_1\cup \Omega_2}\frac{1}{p}|{\nabla u}|^p\,dx\\
			&+\frac{1}{2} \int_{\Omega\setminus \Omega_1\cup \Omega_2}\frac{1}{p}|{\nabla v}|^p - (p_--1)\frac{\varepsilon^2 }{32}\left(\frac{\varrho({ u})+\varrho({ v})}{2}\right)\\ \nonumber &=\frac{1}{2}\left(\varrho({ u})+\varrho({\ v})\right)-(p_--1)\frac{\varepsilon^2 }{32}\left(\frac{\varrho({ u})+\varrho({ v})}{2}\right)\\
			& =  \left(1-(p_--1)\frac{\varepsilon^2}{32}\right)\left(\frac{\varrho({ u})+\varrho({ v})}{2}\right).
		\end{align*}
		Set $\displaystyle \delta(\varepsilon) = \min \left\{\frac{\varepsilon}{2},\ 1-(p_--1)\frac{\varepsilon^2}{32}\right\}$.  Then
		$$\varrho\left(\frac{{ u}+{ v}}{2}\right) \leq \Big(1-\delta(\varepsilon)\Big) \, \frac{\varrho({ u})+\varrho({ v})}{2},$$
		which completes the proof of our claim.
	\end{proof}
	
	\medskip
	We close this section with the following theorem, whose proof can be found in \cite{BMB}.
	
	\medskip
	\begin{theorem}\label{bmb}
		If $1<p_-\leq p<\infty$ in $\Omega$, the modular $\rho_p : L^{p}(\Omega)\rightarrow [0,\infty]$ defined by
		$$\rho_p(u) =\int_{\Omega}\frac{|u|^{p(x)}}{p(x)}\,dx$$
		is $(UC^*)$.
	\end{theorem}
	
	\medskip
	The following corollary follows immediately from Theorems \ref{N=2} and \ref{bmb}:
	
	\begin{corollary}\label{Fisuc}
		If $1<p_-\leq p<\infty$ in $\Omega$, the modular $\rho_{1,p} : W^{1,p}(\Omega)\rightarrow [0,\infty]$ defined by $\rho_{1,p}(u) = \rho_p(u) + \rho_p(\nabla u)$, i.e.,
		$$\rho_{1,p}(u) = \int_{\Omega}\frac{1}{p} |u|^p dx + \int_{\Omega}\frac{1}{p} |\nabla u|^{p}dx $$
		is $(UC^*)$. In particular, $\rho_{1,p}$ is strictly convex, that is, for $u$ and $v$ in $W^{1,p}(\Omega)$, $u\neq v$, one has
		$$\rho_{1,p}\left(\frac{u+v}{2}\right)< \frac{\rho_{1,p}(u)+\rho_{1,p}(v)}{2}.$$
	\end{corollary}
	
	\medskip
	\section{The Dirichlet Integral for unbounded $p$} \label{s6}
	This section focuses on the relationship between modulars and PDEs. Gateaux-like derivatives of modulars will be derived and their connection to the weak formulation of the variable exponent $p$ Laplacian will be demonstrated.  To streamline our proofs, we confine ourselves in this and the following section to variable exponents $p(x)$ that are continuous and unbounded on $\Omega$.\\
	For $\rho_{1,p}$ as defined in (\ref{defrho})
	let $\varphi \in W^{1,p}(\Omega)$ be chosen so that $\rho_{1,p}(\varphi)<\infty$. Consider the functional $F: V^{1,p}_0(\Omega)\rightarrow [0,\infty]$ defined by
	\begin{equation}\label{defF}
		F(u):= \rho_{1,p}(u-\varphi) = \int_{\Omega}\frac{| u-\varphi|^p}{p}\,dx + \int_{\Omega}\frac{|\nabla (u-\varphi)|^p}{p}\,dx.
	\end{equation}
	Notice that $F$ is bounded below and $F(0)<\infty$; more generally, $F\left(C^{\infty}_0(\Omega)\right)\subseteq [0,\infty)$. Indeed, for $u\in C^{\infty}_0(\Omega)$, let $p^*=\sup\limits_{ \supp u} p$. Then $p^*<\infty$ and
	\begin{align}\rho_{1,p}(u-\varphi) &\leq \left(\int\limits_{\text{supp}\,u}+\int\limits_{\Omega\setminus \text{supp}\,u}\right)\frac{1}{p}\left(|\nabla (u-\varphi)|^p+|u-\varphi|^p\right)dx \\ \nonumber & \leq  2^{p^*-1}\int\limits_{\supp u}\frac{1}{p}\left(|\nabla u|^p+ |\nabla \varphi|^p+|u|^p+|\varphi|^p\right)dx \\ \nonumber & + \int\limits_{\Omega \setminus \supp u}(|\nabla \varphi|^p+|\varphi|^p)dx<\infty,
		\end{align}
	  Hence $d:=\inf\limits_{u\in V_0^{\infty}(\Omega)}F(u)<\infty$.  Select a minimizing sequence $(u_j)\subset V^{1,p}_0(\Omega)$, i.e., $\lim\limits_{j \to \infty} F(u_j) = d$.  Without loss of generality, we assume $d > 0$.  We claim that $\displaystyle \left(\frac{u_j}{2}\right)$ $\rho_{1,p}$-converges in $W^{1,p}(\Omega)$.  To see this, observe that if $\displaystyle \left(\frac{u_j}{2}\right)$ were not $\rho_{1,p}-$Cauchy, then for some $\varepsilon>0$ there would be sequences $(i_k)$, $(j_k)$ such that, for all $k\in {\mathbb N},$ it would hold $\displaystyle \rho_{1,p}\left(\frac{u_{j_k}-u_{i_k}}{2}\right) \geq \varepsilon$. By assumption, there exists $L\in {\mathbb N}$ such that $F(u_{j})=\rho_{1,p}(\varphi-u_j)<d(1+\varepsilon)$ for  $j\geq L$.  By virtue of Definition \ref{defuc} and Corollary \ref{Fisuc}, it follows that there exists $\eta \in (0,1)$ such that for $k$ chosen in such away that $j_k,i_k\geq L$, it must hold that (notice that the first inequality follows from the convexity of $V^{1,p}_0(\Omega)$)
	$$d\leq F\left(\frac{u_{j_k}+u_{i_k}}{2}\right) \leq (1-\eta)\, \left(\frac{F(u_{j_k}) + F(u_{i_k})}{2}\right).$$
	If we let $k\rightarrow \infty$, we get $d \leq (1-\eta) d$, which is a contradiction.  Hence $\displaystyle \left(\frac{u_j}{2}\right)$ is $\rho_{1,p}$-Cauchy in $W^{1,p}(\Omega)$.  Since $W^{1,p}(\Omega)$ is $\rho_{1,p}$-complete, then $\displaystyle \left(\frac{u_j}{2}\right)$ $\rho_{1,p}$-converges to $\displaystyle \frac{u}{2} \in W^{1,p}(\Omega)$.  Since $V^{1,p}_0(\Omega)$ is a $\rho_{1,p}$-closed subspace of $W^{1,p}(\Omega)$, we conclude that $u \in V^{1,p}_0(\Omega)$.  We claim that $u$ minimizes $F$, that is, $F(u) = d$.  Indeed, on account of Corollary \ref{aeconvergence1}, since $\displaystyle \frac{u_j}{2}\xrightarrow{\rho_{1,p}}\frac{u}{2}$, one can extract a subsequence $\displaystyle \left(\frac{u_{j_k}}{2}\right)$ that converges to $\displaystyle \frac{u}{2}$ a.e. in $\Omega$, i.e., $u_{j_k}\xrightarrow{\rho_{1,p}}u$ a.e. in $\Omega$.
	Fatou's lemma then yields
	\begin{align*}
		d&\leq \rho_{1,p}(\varphi-u)\leq \liminf\limits_{k}\rho_{1,p}\left(\varphi-\left(\frac{u_{j_k}}{2}+\frac{u}{2}\right)\right)\\  & \leq \liminf\limits_{k}\, \liminf\limits_{l}\rho_{1,p}\left(\varphi -\left(\frac{u_{j_k}}{2}+\frac{u_{j_l}}{2}\right)\right)\\
		&\leq \liminf\limits_{k}\, \liminf\limits_{l} \frac{1}{2}\Big[ \rho_{1,p}\left(\varphi-u_{j_k}\right) +\rho_{1,p}(\varphi-u_{i_k}) \Big] = d,
	\end{align*}
	i.e., $F(u) = d$ as claimed.
	
	\medskip
	\begin{lemma}\label{samelimitsimple}
		There exists a unique minimizer of the functional $F$ in $V^{1,p}_0(\Omega).$
	\end{lemma}
	\begin{proof}
		This follows from the strict convexity of $\rho_{1,p}$.  Indeed, if $u$ and $v$ are minimizers of $F$ in $V^{1,p}_0(\Omega)$, then since $V^{1,p}_0(\Omega)$ is a subspace of $W^{1,p}(\Omega)$, one has
		\begin{align*}
			d &\leq F\left(\frac{u+v}{2}\right) = \rho_{1,p}\left(\varphi- \frac{u+v}{2}\right)\\
			&\leq \frac{\rho_{1,p}(\varphi - u) + \rho_{1,p}(\varphi - v)}{2} \\
			& = \frac{F(u) + F(v)}{2} = d.
		\end{align*}
		On account of the strict convexity of ${\rho}_{1,p}$, the implied equality above is only possible if $u=v$.
	\end{proof}

	\section{The $p$-Laplacian} \label{s7}
	Let $\Omega\subset {\mathbb R}^n$ be bounded and smooth domain ($C^2$ would suffice).  Assume, as in the previous section, $p\in C(\Omega)$ with $p_- > n$.  Obviously, we have the norm inclusion
	\begin{equation*}
		W^{1,p}(\Omega)\hookrightarrow W^{1,p_-}(\Omega)
	\end{equation*}
	and by virtue of Theorem \ref{pqmodularembedding1}, the inclusion $i_{p,p_-}:W^{1,p}(\Omega) \hookrightarrow W^{1,p_-}(\Omega)$ is modularly continuous. Consequently, if $K$ is a norm- closed (if and only if $\rho_{1,p_-}$-closed, since $p_-$ is constant) subset of $W^{1,p_-}(\Omega)$ that contains $C^{\infty}_0(\Omega)$ it follows that $i_{p,p_-}^{-1}(K)$ is $\rho_{1,p}$-closed and $C^{\infty}_0(\Omega) \subseteq i_{pp_-}^{-1}(K)$. By definition, $V^{1,p}_0(\Omega)\subseteq K\cap W^{1,p}(\Omega)\subseteq K$.  Thus $V^{1,p}_0(\Omega)$ is contained in any norm-closed subset of $W^{1,p_-}(\Omega)$ that contains $C^{\infty}_0(\Omega)$. In all,
	\begin{equation}\label{inclusionv1pw1p-}
		V^{1,p}_0(\Omega)\subseteq W^{1,p_-}_0(\Omega).
	\end{equation}
	Consider now a sequence $(v_j)\subset V^{1,p}_0(\Omega)$ such that $(\nabla v_j)$ is $\rho_p$-Cauchy in $\left(L^p(\Omega)\right)^n$. On account of the modular continuity of the inclusion $L^p(\Omega)\subseteq L^{p_-}(\Omega)$,  $(\nabla v_j)$ is $\rho_{p_-}$-Cauchy in $\left(L^{p_-}(\Omega)\right)^n$ and by virtue of (\ref{inclusionv1pw1p-}), the continuity of the Sobolev embedding
	\begin{equation*}
		W^{1,p_-}(\Omega)\hookrightarrow C(\overline{\Omega})
	\end{equation*}
	and the Poincar\'{e} inequality on $W^{1,p_-}_0(\Omega)$, $(v_j)$ is a Cauchy sequence in $W^{1,p_-}_0(\Omega)$ and for any subindices $i,j$, one has (for some constant $C(p_-,\Omega)>0$):
	\begin{equation*}
		\|v_j-v_i\|_{C(\overline{\Omega})} \leq C(p_-,\Omega) \, \| (v_j-v_i)\|_{W^{1,p_-}(\Omega)} \leq C(p_-,\Omega) \, \||\nabla (v_j-v_i)|\|_{L^{p_-}(\Omega)}.
	\end{equation*}
	It follows then that
	\begin{equation*}
		\lim_{i,j \to \infty}\, \rho_p(v_j-v_i) = \lim_{i,j \to \infty}\, \int_{\Omega}\frac{|v_j-v_i|^{p}}{p} \, dx = 0.
	\end{equation*}
	By assumption, one then has
	\begin{equation*}
		\lim_{i,j \to \infty}\, \int_{\Omega} \frac{|v_j-v_i|^p}{p}\, dx + \int_{\Omega} \frac{|\nabla (v_j-v_i)|^p}{p}\, dx =  0,
	\end{equation*}
	and taking into account the fact that $V^{1,p}_0(\Omega)$ is a modularly closed subspace of $W^{1,p}(\Omega)$, it follows that the sequence $(v_j)$ converges modularly in $W^{1,p}(\Omega)$ to a function $v$ in $V^{1,p}_0(\Omega).$  Thus, the following lemma holds:
	
	\medskip
	\begin{lemma}\label{mmain}
		Any sequence $(v_j)\subset V^{1,p}_0(\Omega)$, with $(\nabla v_j)$ $\rho_p$-Cauchy in $\left(L^p(\Omega)\right)^n$, must $\rho_p$-converge in $L^p(\Omega)$ to a function $v\in V^{1,p}_0(\Omega).$ Therefore, $(v_j)$ must converge to $v\in V^{1,p}_0(\Omega)$ in the modular topology of $W^{1,p}(\Omega)$.
	\end{lemma}
	
	\medskip
	We are now ready to prove the following:
	
	\begin{theorem}\label{Dirichletvarphi}
		Let $\Omega \subseteq {\mathbb R}^n$ be a $C^{2}$, bounded domain whose boundary will be denoted by $\partial\Omega$, $p\in C({\Omega})$, $p_->n$ (notice this includes $p_+=\infty$). Then for any $\varphi\in W^{1,p}(\Omega)$ such that  $\int\limits_{\Omega}\frac{|\varphi|^p}{p}dx<\infty$, there exists a solution $u\in W^{1,p}(\Omega)$ to the Dirichlet problem
		\begin{equation}\label{DiPr}
			\begin{cases}
				\Delta_p u=0,\\
				u|_{\partial \Omega}=\varphi.
			\end{cases}
		\end{equation}
		The boundary condition is to be interpreted as $u-\varphi\in V^{1,p}_0(\Omega).$
	\end{theorem}
	\begin{proof}
		Let us consider the functional $\mathcal{F} : V^{1,p}_0(\Omega)\rightarrow [0,\infty]$ defined by
		$$\mathcal{F}(u) = \int_{\Omega}\frac{|\nabla (u-\varphi)|^p}{p}\, dx.$$
		As we did in the previous section (see Lemma \ref{samelimitsimple}), we use the modular uniform convexity $(UC^*)$ to show that $\mathcal{F}$ has a unique minimizer in $V^{1,p}_0(\Omega)$.  Indeed, the functional $F$ is obviously bounded below on $V^{1,p}_0(\Omega)$ and the choice of $\varphi$ guarantees that $\inf\limits_{v\in V_0^{1,p}(\Omega)}F(v)<\infty$. Take a minimizing sequence $(u_j)\subset V^{1,p}_0(\Omega)$, that is, 
		$$\lim_{j \to \infty}\ \mathcal{F}(u_j) = \inf\limits_{u\in V^{1,p}_0(\Omega)}F(u) = d.$$
		Using the uniform convexity $(UC^*)$, we show that $\displaystyle \left(\frac{u_j}{2}\right)$ is $\varrho_p$-Cauchy in $L^p(\Omega)$.  By virtue of Lemma \ref{mmain}, $\displaystyle \left(\frac{ u_j}{2}\right)$ must $\rho_{1,p}$-converge to a function $u\in W^{1,p}(\Omega)$. Since  $V^{1,p}_0(\Omega)$ is closed in the modular topology of $W^{1,p}(\Omega)$, it follows that $u\in V^{1,p}_0(\Omega)$, which is a vector subspace of $W^{1,p}(\Omega)$, which yields $2u\in V^{1,p}_0(\Omega)$.  On account of Lemma \ref{aeconvergence}, it can be assumed, without loss of generality, that  $\displaystyle \left(\frac{\nabla u_j}{2}\right)$ converges pointwise $a.e.$ to $\nabla u$. Fatou's Lemma then yields
		$$\int_{\Omega}\left |\nabla (2u-\varphi)\right|^pdx\leq \liminf_{j \to \infty} \int_{\Omega}\left |\nabla(u_j-\varphi)\right|^p\, dx.$$
		It is claimed that $2u$ is a minimizer of $\mathcal F$ in $K$. As we did in the previous section, we have
		\begin{align*}
			d& \leq \int_{\Omega}\frac{\left|\nabla\left(\varphi-2u\right)\right|^p}{p}dx \leq \liminf\limits_{k\to \infty}\, \int_{\Omega}\frac{\left|\nabla\left(\varphi -\left(\frac{u_{k}}{2}+u\right)\right)\right|^p}{p}dx \\
			& \leq \liminf\limits_{k \to \infty}\, \liminf\limits_{l \to \infty}\, \int_{\Omega}\, \frac{\left|\nabla \left(\varphi-\left(\frac{u_{k}}{2} + \frac{u_{l}}{2}\right)\right)\right|^p}{p}\, dx\\ 
			&\leq \liminf\limits_{k \to \infty}\, \liminf\limits_{l \to \infty}\,  \frac{1}{2}\left(\int_{\Omega}\frac{\left|\nabla\left(\varphi-u_{k}\right)\right|^p}{p}\, dx+ \int_{\Omega}\frac{\left|\nabla \left(\varphi-u_l\right)\right|^p}{p}\, dx\right)\\
			&\leq \liminf\limits_{k \to \infty}\, \liminf\limits_{l \to \infty}\, \frac{1}{2}\, \Big(\mathcal{F}(u_k) + \mathcal{F}(u_l)\Big) = d,
		\end{align*}
		i.e., $2u$ is a minimizer of $\mathcal F$ on $V^{1,p}_0(\Omega).$  Also in this case, two arbitrary minimizing sequences must converge to the same limit. The proof is similar to that of Lemma \ref{samelimitsimple} and will be omitted.  The rest of the proof follows by observing that the $p$-Laplacian $\Delta_p$ is the Gateux-type derivative of $\mathcal F$.
	\end{proof}
	
	\medskip
	\noindent Observe that in particular, if the minimizer $2u$ obtained in the preceding Lemma satisfies $\displaystyle \int_{\Omega}\left|\nabla(2u-\varphi)\right|^p\, dx<\infty$, then it has the following property:
	
	\medskip
	\begin{theorem}
		For any $v\in V^{1,p}_0(\Omega)$ such that $\displaystyle \int_{\Omega}|\nabla (v-\varphi)|^p\, dx <\infty$, it holds
		\begin{equation*}\label{conditionforuniq}
			\int_{\Omega}|\nabla (2u-\varphi)|^{p-2}\ \nabla (2u-\varphi)\ \nabla (v-2u)\, dx \geq 0.
		\end{equation*}
	\end{theorem}
	\begin{proof}
		For any $t \in [0,1]$, the convexity yields
		\begin{align*}
			\int_{\Omega}\frac{1}{p}\left|\nabla\left(2u+t(v-2u)-\varphi\right)\right|^p\, dx& = \int_{\Omega}\frac{1}{p}\left|\nabla \left((1-t)(2u-\varphi)+t(v-\varphi)\right)\right|^p\, dx\\
			&\leq \int_{\Omega}\frac{1-t}{p}\ \left|\nabla(2u-\varphi)\right|^p\, dx + \int_{\Omega}\frac{t}{p}\ \left|\nabla (v-\varphi)\right|^p\, dx\\ 
			& < \infty.
		\end{align*}
		Set
		$$g(t) = \frac{1}{pt}\Big(\left|\nabla\left(2u+t(v-2u)-\varphi\right)\right|^p-|\nabla (2u-\varphi)|^p\Big).$$
		It is straightforward to check that $\lim\limits_{t \to 0^+}\ g(t) = |\nabla(2u-\varphi)|^{p-2}\nabla(2u-\varphi))\nabla (v-2u)$ pointwise.  On the other hand for $t \in [0,1]$, the convexity and the minimal character of $2u$ yield
		\begin{align*}
			0\leq g(t)&=\frac{1}{pt}\Big(\left|\nabla\left((1-t)(2u-\varphi)+t(v-\varphi)\right)\right|^p-|\nabla (2u-\varphi)|^p\Big)\\ 
			&\leq\ \frac{1}{pt}\Big((1-t)|\nabla (2u-\varphi)|^p+t|\nabla (v-\varphi)|^p-|\nabla (2u-\varphi)|^p\Big)\\ 
			&= \frac{1}{p}\Big(|\nabla (v-\varphi)|^p-|\nabla (2u-\varphi)|^p\Big),
		\end{align*}
		which, by assumption, is integrable. As an immediate consequence, one has
		\begin{equation*}
			\lim_{t\rightarrow 0^+}\ \frac{\mathcal{F}\left(2u+t(v-2u)\right)-\mathcal{F}(2u)}{t} = \int_{\Omega}|\nabla (2u-\varphi)|^{p-2}\nabla(2u-\varphi)\nabla(v-2u)\, dx.
		\end{equation*}
		The integral in the right-hand side is finite and since $2u$ is assumed to be a minimizer of $\mathcal{F}$ on $V^{1,p}_0(\Omega)$, the last equality forces
		\begin{equation*}\label{condforuniq}
			\int_{\Omega}|\nabla (2u-\varphi)|^{p-2}\ \nabla(2u-\varphi)\ \nabla \left(v-2u\right)\, dx \geq 0,
		\end{equation*}
		as claimed.
	\end{proof}
	
	\medskip
	We underline the following consequence which is of independent interest
	
	\begin{corollary} Under the hypothesis of Theorem \ref{Dirichletvarphi}, consider the operator $S: V^{1,p}_0(\Omega) \rightarrow \left(C^{\infty}_0(\Omega)\right)^{\ast}$ defined by
		\begin{equation*}
			\langle S(w), h\rangle=\int_{\Omega}|\nabla (w-\varphi)|^{p-2}\ \nabla (w-\varphi)\ \nabla h\, dx.
		\end{equation*}
		Then there exists a solution to the variational inequality
		$$\langle S(w), v-w\rangle\ \geq\ 0\,$$
		for all $v\in V^{1,p}_0(\Omega)$, with $\displaystyle \int_{\Omega}\ \frac{|\nabla (v-\varphi)|^p}{p}\, dx < \infty$ satisfying $\displaystyle \int_{\Omega}\ \frac{|\nabla (u-\varphi)|^p}{p}\, dx < \infty.$
	\end{corollary}
	
	\medskip
	Next we set out to prove that there is a unique weak solution $v\in V^{1,p}_0(\Omega)$ of problem (\ref{DiPr}). Specifically:
	
	\begin{theorem}\label{uniqueness}
		There exists at most one weak solution $v\in W^{1,p}(\Omega)$ to the Dirichlet problem (\ref{DiPr}) satisfying $\displaystyle \int_{\Omega}p^{-1}\ |\nabla v|^p\, dx < \infty$ and such that the inequality
		\begin{equation}\label{additionalcondition}
			\int_{\Omega}\ |\nabla v|^{p-2}\ \nabla v\ \nabla (w+v-\varphi)\, dx \leq 0
		\end{equation}
		holds for every $w\in V^{1,p}_0(\Omega)$ such that $\displaystyle \int_{\Omega}p^{-1}\ |\nabla (w-\varphi)|^p\, dx < \infty.$
	\end{theorem}
	\begin{proof}
		To see this, let $v_1, v_2$ be two such weak solutions and recall the inequality (valid for all vectors $A\in {\mathbb R}^n, B\in {\mathbb R}^n$, $p\geq 2$ and a positive constant $\gamma(p)$ \cite{Lindqvist})
		\begin{equation*}
			\langle |A|^{p-2}A-|B|^{p-2}B,A-B\rangle\geq \gamma(p)|A-B|^p,
		\end{equation*}
		where the brackets represent the standard inner product in ${\mathbb R}^n$.
		Then $v_1$ and $v_2$ satisfy inequality (\ref{additionalcondition}) with $w_1= \varphi-v_1$ and $w_2=\varphi-v_2$, respectively, that is
		\begin{equation*}
			\int_{\Omega}|\nabla v_1|^{p-2}\nabla v_1\nabla (v_1-v_2)\, dx\leq 0\;\;\; and\;\;\; \int_{\Omega}|\nabla v_2|^{p-2}\nabla v_2\nabla (v_2-v_1)\, dx\leq 0.
		\end{equation*} 
		Hence
		\begin{align*}
			\int_{\Omega}\gamma(p)|v_1-v_2|^p \, dx & \leq \int_{\Omega}\Big(|\nabla v_1|^{p-2}\nabla v_1\nabla (v_1-v_2)-|\nabla v_2|^{p-2}\nabla v_2\nabla (v_1-v_2)\Big)\, dx\\ 
			&= \int_{\Omega}\left(|\nabla v_1|^{p-2}\nabla v_1-|\nabla v_2|^{p-2}\nabla v_2\right)\nabla (v_1-v_2)\, dx \\
			&\leq 0,
		\end{align*}
		from which the claim follows at once.
	\end{proof}
	
	\medskip
	This section is closed with the following theorem:
	
	\begin{theorem}\label{Dirichlet+}
		Let $\Omega\subset {\mathbb R}^n$ be a bounded domain with  smooth boundary $\partial\Omega$, $p\in C(\Omega)$, $p_->n$ (notice that the possibility $p_+=\infty$ is included here). Then, for any $\varphi \in W^{1,p}(\Omega)$ with $\|\varphi\|_{1,p}\leq 1$, any nonegative $q\in C(\overline{\Omega})$, there exists a unique solution $u\in W^{1,p}(\Omega)$ to the boundary value problem
		\begin{equation*}\label{general}
			\begin{cases}
				-\Delta_p(u)+q\ |u|^{p-2}u = 0 \;\;\text{in}\;\;\; \Omega\\
				u|_{\partial \Omega}= \varphi,
			\end{cases}
		\end{equation*}
		such that $\displaystyle \int_{\Omega}\Big(\left|\nabla u\right|^p+|u|^p\Big)\, dx < \infty$ and 
		$$\int_{\Omega}\Big(|\nabla v|^{p-2}\ \nabla v\ \nabla (w+v-\varphi)+ |v|^{p-2}\ v\ (w+v-\varphi)\Big)\, dx \leq 0,$$
		for each $w\in V^{1,p}_0(\Omega)$ such that $\displaystyle \int_{\Omega}\Big(|\nabla (w-\varphi)|^p + |w-\varphi |^p\Big)\, dx < \infty.$
	\end{theorem}
	\begin{proof}
		The existence proof is similar to that of Theorem \ref{Dirichletvarphi}.  In this case, a solution is given by the unique minimizer of the functional $J : V^{1,p}_0(\Omega)\rightarrow [0,\infty]$ defined by
		$$ J(u) = \int_{\Omega}\left(\frac{1}{p}|\nabla (u-\varphi)|^p+\frac{q}{p}|u-\varphi|^p\right)\, dx.$$
		The existence and uniqueness of the minimizer can be proved by invoking the proper version of Corollary \ref{Fisuc} (here the condition $q\geq0$ is essential) along the same lines as those of the proof of Lemma \ref{samelimitsimple}. The uniqueness part follows as in Theorem \ref{uniqueness}.
	\end{proof}
	
	\medskip
	\section{A remark on the connection with the classical $p$-Laplacian}\label{Sectionremark}
	If $p$ is constant on $\Omega$, there is no essential distinction between the operators
	\begin{equation}\label{classical}
		\text{div}\left(|\nabla u|^{p-2}\nabla u\right)
	\end{equation}
	and
	\begin{equation}\label{pclassical}
		\text{div}\left(p|\nabla u|^{p-2}\nabla u\right).
	\end{equation}
	For variable $p$, though, (\ref{classical}) and (\ref{pclassical}) are very different operators and both generalize the $p$-Laplacian. Therefore, it is in order to remark that the Dirichlet problem for the operator (\ref{pclassical}) can be handled along lines similar to those depicted above.\\
	Specifically, in order to deal with the Dirichlet problem
	\begin{equation}\label{DPP}
		\begin{cases}
			\text{div}\left(p|\nabla u|^{p-2}\nabla u\right)=0 \;\;\; \text{in} \;\;\; \Omega\\
			u|_{\partial \Omega}= \varphi
		\end{cases}
	\end{equation}
	the modular defined in (\ref{def-rho}) has to be replaced with $\displaystyle \rho(u) = \int_{\Omega}|u|^p \, dx$. The results in Section \ref{s4} carry over with no difficulties to the new modular. In fact, Lemma \ref{pimpliesq} holds in this case without the additional condition (\ref{condpq}). The spaces $V^{1,p}_0(\Omega)$ and $U^{1,p}_0(\Omega)$ can be defined in terms of $\rho$. Section \ref{s5} carries over to the modular $\rho$ almost directly. The solution to the problem (\ref{DPP}) entails the minimization of the Dirichlet integral $G : V^{1,p}_0(\Omega)\rightarrow [0,\infty]$ where
	\begin{equation}
		G(u) = \int_{\Omega}|\nabla (u-\varphi)|^p\ dx,
	\end{equation}
	which can be achieved using {\it mutatis mutandis} along the same lines as those in Section \ref{s5}. These considerations lead to the following result:
	\begin{theorem}
		Let $\Omega \subset {\mathbb R}^n$ be a bounded domain with a smooth boundary $\partial\Omega$, $p\in C(\Omega)$ with $p_-=\inf\limits_{x \in \Omega}\, p(x) > n$ and $\varphi \in W^{1,p}(\Omega)$ such that $\displaystyle \int_{\Omega}|\nabla \varphi|^p \, dx< \infty.$
		Then there exists a unique solution $u\in W^{1,p}(\Omega)$ to the Dirichlet problem (\ref{DPP}) that satisfies
		\begin{equation}\label{additionalconditionp}
			\int_{\Omega}|\nabla u|^{p-2}\nabla u\nabla (w+u-\varphi)dx\leq 0
		\end{equation}
		for every $w\in V^{1,p}_0(\Omega)$ such that $\displaystyle \int_{\Omega}|\nabla (w-\varphi)|^p \, dx < \infty.$
	\end{theorem}
	
	\medskip
	Problem (\ref{general}) can also be reformulated and solved for this version of the $p$-Laplacian in the same manner. The details will be omitted.
	
	\medskip
	\section*{Acknwoledgment}
	The third author would like to express his appreciation to the Department of Mathematics at Khalifa University, where part of this investigation was carried out.

\end{document}